\documentclass[11pt]{amsart}

\usepackage[T1]{fontenc}
\usepackage[utf8]{inputenc}
\usepackage{lmodern}
\usepackage[margin=1.08in]{geometry}
\usepackage{amsmath,amssymb,amsthm,mathtools,mathrsfs}
\usepackage{microtype}
\usepackage{xcolor}
\usepackage[colorlinks=true,linkcolor=blue,citecolor=green!45!black,urlcolor=blue]{hyperref}
\hypersetup{
  pdftitle={Generalized Whittaker models and rational orbit closures for reductive groups},
  pdfauthor={Taiwang Deng},
  pdfsubject={An extension of Gomez--Gourevitch--Sahi Theorem E to connected reductive groups},
  pdfkeywords={generalized Whittaker models, rational nilpotent orbits, local character expansion, reductive groups}
}

\usepackage[english]{babel}
\usepackage{stackengine}
\usepackage{scalerel}

\allowdisplaybreaks
\numberwithin{equation}{section}
\theoremstyle{plain}
\newtheorem{theorem}{Theorem}[section]
\newtheorem{proposition}[theorem]{Proposition}
\newtheorem{lemma}[theorem]{Lemma}
\newtheorem{corollary}[theorem]{Corollary}

\theoremstyle{definition}

\theoremstyle{remark}
\newtheorem{remark}[theorem]{Remark}

\newcommand{\C}{\mathbb C}
\newcommand{\g}{\mathfrak g}
\newcommand{\SL}{\operatorname{SL}}

\newcommand{\GL}{\operatorname{GL}}
\newcommand{\Ad}{\operatorname{Ad}}
\newcommand{\ad}{\operatorname{ad}}
\newcommand{\Hom}{\operatorname{Hom}}
\newcommand{\End}{\operatorname{End}}
\newcommand{\Wh}{\operatorname{Wh}}
\newcommand{\Tr}{\operatorname{Tr}}
\newcommand{\vol}{\operatorname{vol}}
\newcommand{\supp}{\operatorname{supp}}
\newcommand{\OO}{\mathcal O}
\newcommand{\one}{\mathbf 1}
\newcommand{\clan}[1]{\overline{#1}}

\title[Whittaker models and rational orbit closures]{Generalized Whittaker models and rational orbit closures for reductive groups}
\author{Taiwang Deng}
\address[T.D.]{Beijing Institute of Mathematical Sciences and Applications (BIMSA), Huairou District, 100084, Beijing\\
China}
\email{dengtaiw@bimsa.cn}
\subjclass[2020]{Primary 22E50; Secondary 11F70, 20G25}
\keywords{Generalized Whittaker models, rational nilpotent orbits, local character expansions, reductive groups, Slodowy slices}
\date{}

\begin{document}

\begin{abstract}
Let $F$ be a finite extension of $\mathbb Q_p$, let $G$ be a
connected reductive group over $F$, and let $\pi$ be an admissible
smooth complex representation of $G(F)$ of finite length.
We prove that a rational nilpotent orbit has a nonzero
neutral generalized Whittaker model on $\pi$ if and only if the orbit
lies in the union of the $F$-analytic closures of the orbits with nonzero
coefficient in the local character expansion of $\pi$. In the case of
$G=\GL_n$, this equivalence was proved by G\'omez, Gourevitch,
and Sahi.

For the converse, vanishing of the corresponding generalized
Whittaker space yields a test function that acts by zero on $\pi$,
but whose trace is nonzero by the local character expansion,
giving a contradiction. The construction combines rational
Slodowy slices with Fourier analysis on small compact open subgroups.
\end{abstract}

\maketitle
\tableofcontents

\section{Introduction}\label{sec:introduction}

\subsection{Background and main result}

Let $F/\mathbb Q_p$ be finite, and fix a nontrivial continuous
unitary additive character $\psi:F\to\C^\times$.
Let $G$ be a connected reductive algebraic group over $F$.
We use the same letter for the group of rational points,
$G=G(F)$, and put $\g=\operatorname{Lie}(G)(F)$.
Throughout, a rational orbit means a $G(F)$-orbit in $\g$.
We take closures in the $F$-analytic topology of $\g$ and
denote the closure of $\OO$ by $\clan{\OO}$.

The local character expansion (see Section \ref{subsec: local-character}) and generalized Whittaker models
attach nilpotent orbits to an admissible smooth representation
of $G$ of finite length.
The character expansion records Fourier transforms of measures supported on nilpotent
orbital near the identity and Whittaker models record
equivariant functionals under unipotent subgroups. The relation at
maximal orbits is due to M{\oe}glin and Waldspurger \cite{MW87},
with the extension to residue characteristic two by Varma \cite{Varma14}.
We follow the formulation in \cite[Theorem~4.3.1]{GGS17}.

For general linear groups, G\'omez, Gourevitch, and Sahi prove
the complete relation in \cite[Theorem~E]{GGS17}: a neutral
generalized Whittaker model is nonzero precisely when the
associated orbit lies in the closure of the union of the orbits
with nonzero coefficients in the local character expansion.
We prove the same relation for every
connected reductive $F$-group. The distinction between rational
orbits is essential: a geometric nilpotent orbit, meaning an orbit
over an algebraic closure of $F$, may contain several
$G(F)$-orbits with different analytic closures
\cite[Section~6 and Lemma~44]{McNinch}.
In this paper we present a generalization of the result of G\'omez, Gourevitch, and Sahi 
for general reductive group over $F$ which we state below.

Fix a nondegenerate $G$-invariant symmetric $F$-bilinear form $B$
on $\g$ which identifies adjoint and
coadjoint orbits. For an admissible smooth complex representation $\pi$ of finite length,
let $c_{\OO}(\pi)$ be the coefficient attached to $\OO$ in the
local character expansion, with the conventions of
Section~\ref{sec:setup}. Following \cite[Section~1.1]{GGS17},
write $[\OO]$ for the formal symbol indexed by $\OO$ and define
the wave-front cycle and wave-front set by
\[
 \operatorname{WFC}(\pi)=\sum_{\OO}c_{\OO}(\pi)[\OO],
 \qquad
 \operatorname{WF}(\pi)=
 \bigcup_{c_{\OO'}(\pi)\ne0}\clan{\OO'}.
\]
The sum and union range over rational nilpotent orbits. Note that there
are finitely many rational nilpotent orbits in $ \operatorname{WF}(\pi)$ by
\cite[Theorem~40]{McNinch}.

For a rational nilpotent orbit $\OO$, write $\Wh_{\OO}(\pi)$ for
the associated neutral generalized Whittaker space.
When a representative $f\in\OO$ is specified, we also write
$\Wh_f(\pi)$. Section~\ref{subsec:neutral-model} gives the precise
construction, including the convention for the zero orbit.

\begin{theorem}\label{thm:main}
For every admissible smooth complex representation $\pi$ of $G$
of finite length and every rational nilpotent orbit $\OO$,
\begin{equation}\label{eq:main}
 \Wh_{\OO}(\pi)\ne0
 \quad\Longleftrightarrow\quad
 \OO\subset\operatorname{WF}(\pi).
\end{equation}
\end{theorem}

Our motivation comes in part from our study of degenerate Whittaker
models for $\operatorname{GL}_{2n}(F)$ in \cite{Deng26}.
The criterion in \cite[Theorem~E]{GGS17} controls the vanishing of
the twisted Jacquet module for the parabolic of type $(n,n)$, associated with the
nilpotent orbit of Jordan type $(2^n)$. The orbit criterion can then
be expressed in Langlands--Zelevinsky multisegment data
\cite[Lemma~1.3 and Theorems~1.4 and~2.5]{Deng26}.
Theorem~\ref{thm:main} extends the underlying orbit closure
criterion beyond general linear groups.

Related work of Jiang, Liu, and Zhang \cite{JLZ25} studies rational
wave-front sets for classical groups through the local Langlands
correspondence. A local Langlands parameter is called \emph{generic}
if its packet on the quasi-split group contains a representation
with a nonzero Whittaker model attached to a regular nilpotent orbit.
For irreducible admissible representations of quasi-split classical
groups and their pure inner forms with such parameters,
Jiang, Liu, and Zhang construct an arithmetic wave-front set by successive descents of
the parameter together with the additional data specifying the
representation. Jiang, Liu, Luo, Ma, and Zhang prove that every
orbit in this arithmetic wave-front set has a nonzero generalized
Whittaker model \cite[Theorem~1.2 and Corollary~3.6]{JLLMZ26}.
The character comparisons in these works concern maximal orbits.
In particular, \cite[Theorem~1.5]{JLLMZ26} treats real classical
groups. 

Theorem~\ref{thm:main} should be seen as a special case of 
Liu and Prasad's \cite[Conjecture~1.1]{LP26}.

\subsection{The argument}

The
implication from nonzero Whittaker models to containment in the
wave-front set is known combining work of Moeglin-Walspurger and Varma, 
cf. \cite[Theorem~4.3.1(i)]{GGS17}. Our main result concerns the converse direction which we now explain. Choose a maximal
orbit $\OO'$ occurring in $\operatorname{WFC}(\pi)$ whose closure
contains a nonzero $f$. Choose a rational $\mathfrak{sl}_2$-triple
$(e,h,f)$, with $[h,e]=2e$, $[h,f]=-2f$ and $[e,f]=h$
\cite[Chapter~VIII, \S11, no.~2, Proposition~2]{Bourbaki789}.
Write $\g^x=\ker(\ad x)$ for the
Lie-algebra centralizer of $x\in\g$. The Slodowy slice
$S_f=f+\g^e$, formed using the chosen triple, meets $\OO'$ by
Lemma~\ref{lem:slice-incidence}. Choose $s\in\OO'\cap S_f$.

Assume that $\Wh_f(\pi)=0$. This assumption first yields a compactly
supported measure $\phi$ whose integrated action annihilates a prescribed
finite-dimensional subspace of $\pi$. Here $\phi$ denotes the measure
$\phi_{C_t,K}$ constructed in Section~\ref{sec:positive}.

For a suitable compact open subset $\mathcal T\subset\g$,
Section~\ref{sec:compact} constructs a function $e_{\mathcal T}$ supported
in a fixed small compact open subgroup $Q$ of $G$ by composing the
inverse Fourier transform of $\one_{\mathcal T}$ with the logarithm.
The results
of Boyarchenko and Sabitova \cite{BS} show that $e_{\mathcal T}$ is an
idempotent (Proposition~\ref{prop:compact-projector}).
By admissibility, $\pi(e_{\mathcal T})$ has
finite-dimensional image, so the preceding measure $\phi$ can be
chosen to satisfy $\pi(\phi*e_{\mathcal T})=0$, where $*$ denotes
convolution. Let $\phi|_Q$ denote the restriction of the measure
$\phi$ to $Q$, we then show
\[
 e_{\mathcal T}*\phi*e_{\mathcal T}
 =e_{\mathcal T}*(\phi|_Q)*e_{\mathcal T}
\]
by Lemma~\ref{lem:localized-sandwich}. All three factors on the
right are supported in $Q$. Thus $e_{\mathcal T}*\phi*e_{\mathcal T}$ is
supported in $Q$ and acts by zero on $\pi$. Averaging this
function over conjugation by $Q$ gives a function with the same
two properties. After pullback by the exponential map, the Fourier
transform is nonnegative and has a positive orbital integral on
$\OO'$. Maximality of $\OO'$ allows the support of this Fourier
transform to avoid every other orbit in $\operatorname{WF}(\pi)$. The character expansion therefore gives a nonzero
trace for the zero operator, a contradiction. The 
convolution formula in \cite[Proposition~2.7(a)]{BS} also permits
this construction in residue characteristic two.

\subsection{Comparison with the proof for general linear groups}
\label{subsec:comparison}

 Our approaches above to prove the converse direction is different from \cite[Theorem~E]{GGS17},  in which the case $G=\GL_n$ is proved.
In the following we compare the two approaches.

A Whittaker pair $(S,s)$ consists of a semisimple element $S\in \g$ with
rational adjoint eigenvalues and a nilpotent element $s\in \g$ satisfying
$[S,s]=-2s$. Write $G_S=Z_G(S)$ for the centralizer of $S$ in $G$.
Denote by $\pi^\vee$ the smooth dual of $\pi$.
The universal degenerate model $\mathcal W_{S,s}$ is the smooth
compactly induced representation attached to the pair in
\cite[Section~2.5]{GGS17}.
The degenerate Whittaker space of $\pi$ associated to $(S,s)$ is
the space of $G$-equivariant linear maps
\[
 \Hom_G(\mathcal W_{S,s},\pi^\vee).
\]
Since $\mathcal W_{S,s}$ is smooth, every $G$-equivariant map
from $\mathcal W_{S,s}$ to $\Hom_{\C}(\pi,\C)$ has image in
$\pi^\vee$: an open subgroup fixing a vector also fixes its image.
The pair $(S,s)$ is neutral if $(e,S,s)$ is a rational
$\mathfrak{sl}_2$-triple for some $e\in\g$.
For a neutral pair, we also use the noation $\mathcal W_{\OO}$, where $\OO$ is
the orbit of $s$.

For $G=\operatorname{GL}_n(F)$, choose a maximal orbit $\OO'$
occurring in $\operatorname{WFC}(\pi)$ with
$\OO\subset\clan{\OO'}$. By \cite[Theorem~D]{GGS17}, there is a
Whittaker pair $(S,s)$ such that $s\in\OO'$ and
$\OO\cap\clan{G_Ss}\ne\varnothing$.
The epimorphism theorem \cite[Theorem~A]{GGS17} then gives
\[
 \mathcal W_{\OO}\twoheadrightarrow\mathcal W_{S,s}.
\]
By \cite[Theorem~4.3.1(ii)]{GGS17}, the degenerate Whittaker space
attached to $(S,s)$ is nonzero. Applying $\Hom_G(-,\pi^\vee)$ embeds
that space into $\Wh_{\OO}(\pi)= \Hom_G(\mathcal W_{\OO},\pi^\vee)$ and proves the converse.

The construction in \cite[Section~4.3]{GGS17} requires a point of
$\OO$ to lie in the closure of the $G_S$-orbit of $s$. This condition
is stronger than containment in the closure of a $G$-orbit.
Failure of the condition for general groups is
the obstacle identified in \cite[Remarks~1.2.1--1.2.2]{GGS17}.
Our proof uses only the rational orbit-closure relation to obtain
the intersection with the Slodowy slice. Vanishing of the
Whittaker space then gives a function whose integrated
action is zero. In the local character expansion for this function,
only one orbit contributes, and its orbital integral is positive.
The resulting trace is nonzero; see Proposition~\ref{prop:positive-test} and
Lemma~\ref{lem:orbit-isolation}.
Here \cite{GGS17} use maximality to obtain a
nonzero degenerate model \cite[Theorem~4.3.1(ii)]{GGS17},
whereas we use maximality to make only one nonzero
coefficient contribute to the local character expansion.

For $G=\operatorname{GL}_n$ and $F/\mathbb Q_p$ finite,
Theorem~\ref{thm:main} recovers \cite[Theorem~E]{GGS17}.
In fact, a smooth representation of
$\operatorname{GL}_n(F)$ has finite length if and only if it is
admissible and finitely generated \cite[Theorem~4.1]{BZ76}.
The theorem of G\'omez, Gourevitch, and Sahi also covers
$F=\mathbb R$ and $F=\mathbb C$.

\subsection{Relation to the work of Barbasch and Moy}
\label{subsec:barbasch-moy}

Barbasch and Moy \cite{BM97} construct test functions on $p$-adic
groups by inflating generalized Gelfand--Graev characters from
finite reductive quotients of parahoric subgroups. These characters
are attached to nilpotent orbits. The formula of Lusztig recalled in
\cite[Proposition~2.6]{BM97} computes their Fourier transforms by
counting group elements carrying a given point into a Slodowy
slice, and therefore gives nonnegative values. The orbit-closure
properties in \cite[Proposition~2.7 and Theorem~4.5]{BM97} yield
a triangular system for the coefficients in the local character
expansion \cite[Section~4.6 and Corollary~4.7]{BM97}.
A smooth representation is said to have depth zero if it has a nonzero
vector fixed by the pro-unipotent radical of some parahoric subgroup.
The calculation applies to depth-zero representations under the
group and residue-characteristic hypotheses of
\cite[Section~4.4]{BM97}.

The positivity argument in Lemma~\ref{lem:measure-positivity}
uses a related averaging procedure: the Fourier transform in
\eqref{eq:subgroup-Fourier} is an indicator function, and averaging
under conjugation preserves nonnegativity. Our proof also requires
a test function whose integrated action on $\pi$ vanishes.
Lemma~\ref{lem:localized-sandwich} preserves this vanishing when
the convolution is restricted to a fixed small compact subgroup.
Convolution with the central idempotent then restricts the Fourier
transform to a neighborhood of a chosen maximal orbit, so that
only one orbit contributes to the character expansion.
The idempotence and convolution identities used in this step
follow from \cite{BS}, as explained in Section~\ref{sec:compact}.

\subsection{Organization}

Section~\ref{sec:setup} constructs the neutral Whittaker model and
fixes the Fourier and measure conventions.
Section~\ref{sec:neutral} constructs a compactly supported measure
from vanishing of the Whittaker space, and
Section~\ref{sec:geometry} proves the required results on Slodowy
slices and conjugation by unipotent elements.
Section~\ref{sec:compact} constructs central idempotents on small
compact open subgroups and proves the Fourier transform identity
for their convolution with conjugation-invariant measures.
Section~\ref{sec:positive}
constructs the function with nonnegative Fourier transform, and
Section~\ref{sec:conclusion} reduces the character expansion to a
single orbital integral and proves the theorem.

\subsection*{Acknowledgements}
\addtocontents{toc}{\protect\setcounter{tocdepth}{1}}
The author benefited greatly from the lectures by Professor D.~Prasad 
at BIMSA in the conference "Automorphic representations, L-functions and Harmonic analysis"  and
the IASM workshop "Wavefront Sets of Representations and Endoscopic Classification".
He also thanks Professor Prasad for reading part of the paper and for suggesting 
improvements.

 The author is supported by the National Natural Science Foundation
of China, No.~12401013 and Beijing Natural Science Foundation, No. 1244042.

The author used ChatGPT (OpenAI) for editorial assistance in improving the language, exposition, and notational inconsistencies of the manuscript. All mathematical content and references were verified by the author, who take full responsibility for the manuscript.

\section{Models, Fourier transforms and character expansions}
\label{sec:setup}

\subsection{The neutral model}
\label{subsec:neutral-model}

In characteristic zero the reductive Lie algebra decomposes as
$\g=\mathfrak z\oplus[\g,\g]$, with $\mathfrak z$ its center and
$[\g,\g]$ semisimple
\cite[Proposition~21.61]{MilneAG},
\cite[Chapter~I, \S6, no.~4, Proposition~5 and the following
remark]{Bourbaki123}.
Fix an invariant form $B$ on $\g$ as follows: take the direct sum of the
Killing form on $[\g,\g]$ and
any nondegenerate symmetric form on $\mathfrak z$. The Killing form is nondegenerate by Cartan's criterion
\cite[Chapter~I, \S6, no.~1, Theorem~1]{Bourbaki123}, and the connected group
$G$ acts trivially on the second summand.

For a nonzero rational nilpotent $f\in [\g,\g]$, choose a rational
$\mathfrak{sl}_2$-triple $(e,h,f)$ as above. Such a triple exists by
the Jacobson--Morozov theorem, applied to the semisimple Lie algebra $[\g,\g]$
\cite[Chapter~VIII, \S11, no.~2, Proposition~2]{Bourbaki789}.
The corresponding Lie-algebra homomorphism integrates to an
$F$-homomorphism $\vartheta:\SL_2\to G$ by
\cite[Theorem~23.70(c)--(e)]{MilneAG}. Set
\[
\gamma_h(u):=\vartheta(\operatorname{diag}(u,u^{-1})).
\]
The map $\gamma_h$ is an $F$-rational cocharacter with differential $h$.
The associated neutral grading is
$\g_j=\{Y\in\g:[h,Y]=jY\}$ for $j\in\mathbb Z$. Write
\[
 \g_{\ge i}=\bigoplus_{j\ge i}\g_j,\qquad
 N_i=\exp(\g_{\ge i})\quad(i\ge1).
\]
We also write $\g_{>0}=\bigoplus_{j>0}\g_j$ and use $Y_j$ for
the component of $Y\in\g$ in $\g_j$.
On each $\g_{\ge i}$, the exponential map is defined everywhere
and is an isomorphism of algebraic varieties onto $N_i$.
The exponential map commutes with conjugation and carries the finite
Baker--Campbell--Hausdorff product to multiplication in $N_i$
\cite[14.31--14.32 and~14.35--14.37]{MilneAG}.

We recall the relevant facts from $\mathfrak{sl}_2$-representation
theory. Let $\mathfrak s=Fe\oplus Fh\oplus Ff\simeq\mathfrak{sl}_2(F)$.
For a finite-dimensional $\mathfrak s$-module $V$, let
\[
V_j=\{v\in V:h\cdot v=jv\}.
\]
By complete reducibility and the classification of irreducible
$\mathfrak{sl}_2(F)$-modules
\cite[Chapter~I, \S6, no.~2, Theorem~2]{Bourbaki123},
\cite[Chapter~VIII, \S1, no.~2, Proposition~2 and its
corollary]{Bourbaki789},
$V$ is a direct sum of irreducible modules. An irreducible module
of highest weight $d\in\mathbb Z_{\ge0}$ has one-dimensional weight
spaces of weights $d,d-2,\ldots,-d$. The action of $f$ is nonzero
from each weight space other than the lowest to the next lower
weight space. Consequently, the action of $f$ gives
\[
 \begin{aligned}
 V_j&\longrightarrow V_{j-2} &&\text{injectively for }j>0,\\
 V_0&\longrightarrow V_{-2} &&\text{surjectively},\\
 V_1&\longrightarrow V_{-1} &&\text{isomorphically}.
 \end{aligned}
\]
For the adjoint $\mathfrak s$-module $V=\g$, these maps are
the restrictions of $\ad f$. 

Define an alternating form on $\g_1$ by
$\omega_f(X,Y)=B(f,[X,Y])$. The form $\omega_f$ is nondegenerate. Indeed,
$\ad f:\g_1\to\g_{-1}$ is an isomorphism, and $B$ pairs these
weight spaces nondegenerately.
Choose a Lagrangian (maximal isotropic) subspace $\ell\subset\g_1$
for $\omega_f$. The relation $[\g_i,\g_j]\subset\g_{i+j}$ shows
that $\ell\oplus\g_{\ge2}$ is a nilpotent Lie subalgebra. Define
\[
 L_f=\exp\!\left(\ell+\bigoplus_{j\ge2}\g_j\right),
 \qquad
 \theta_f(u)=\psi(-B(f,\log u))\quad(u\in L_f).
\]

\begin{lemma}\label{lem:whittaker-character}
The function $\theta_f:L_f\to\C^\times$ is a group character.
\end{lemma}

\begin{proof}
For $X,Y\in\ell$, isotropy of $\ell$ means that
\[
 B(f,[X,Y])=\omega_f(X,Y)=0.
\]
If $X\in\g_i$ and $Y\in\g_j$ with $i,j\ge1$ and $i+j\ge3$,
then $[X,Y]\in\g_{i+j}$ and $B(f,[X,Y])=0$. Since $\log L_f=\ell\oplus\g_{\ge2}$,
bilinearity therefore gives $B(f,[X,Y])=0$ for all $X,Y\in\log L_f$.
The Baker--Campbell--Hausdorff formula
\cite[14.35 and Proposition~14.36]{MilneAG} now yields
\[
 B(f,\log(uv))=B(f,\log u)+B(f,\log v)
 \qquad(u,v\in L_f).
\]
Thus $\theta_f(uv)=\theta_f(u)\theta_f(v)$, so $\theta_f$ is a group character.
\end{proof}

Define the neutral generalized Whittaker  space by
\[
\Wh_f(\pi)=\Hom_{L_f}(\pi,\theta_f).
\]
We write
$\Wh_{\OO}(\pi)=\Wh_f(\pi)$ for the chosen representative $f\in\OO$.
Nonvanishing is independent of the triple, Lagrangian and
representative \cite[Sections~2.3--2.5]{GGS17}.
For the zero orbit, set
$\Wh_{\{0\}}(\pi)=\Wh_0(\pi)=\Hom_{\C}(\pi,\C)$.

For a smooth representation $W$ of a group $H$ and a character
$\theta$ of $H$, write $u\cdot w$ for the action of $u\in H$ on
$w\in W$. Define the twisted coinvariant quotient 
\[
 W_{H,\theta}=
 W/\langle u\cdot w-\theta(u)w:u\in H,\ w\in W\rangle.
\]
Ordinary coinvariants $W_H$ are the case $\theta=1$.
The realization of the neutral generalized Whittaker functional
space using the Lagrangian $\ell$ \cite[Lemma~2.5.2]{GGS17} gives
\begin{equation}\label{eq:Wh-definition}
 \Wh_f(\pi)=(\pi_{L_f,\theta_f})^*.
\end{equation}
Here the superscript $*$ denotes the full algebraic dual.

\subsection{Fourier and measure conventions}

Let $\mathcal O_F$ be the valuation ring of $F$ with a
uniformizer $\varpi$. 
\begin{itemize}
\item A lattice in a finite-dimensional $F$-space
is a finitely generated $\mathcal O_F$-submodule spanning that space.
\item
A lattice $M$ in an $F$-Lie algebra is called a \emph{Lie lattice}
if $[X,Y]\in M$ for every $X,Y\in M$, or equivalently
$[M,M]\subset M$.
\end{itemize}

For an $F$-analytic space $X$, let $C_c^\infty(X)$ be the space of 
compactly supported locally constant complex functions on $X$.
Choose an additive Haar measure $dY$ on $\g$ and define the
Fourier transform of $a\in C_c^\infty(\g)$ by
\begin{equation}\label{eq:Fourier}
 \widehat a(\xi)=
 \int_{\g}a(Y)\psi(-B(\xi,Y))\,dY.
\end{equation}
Normalize the Haar measure $d\xi$ on $\g$ so that Fourier inversion
takes the form
\[
 a(Y)=\int_\g\widehat a(\xi)\psi(B(\xi,Y))\,d\xi.
\]
For a finite compactly supported complex measure $\nu$ on $\g$,
define $\widehat\nu$ by the same formula, with integration against
$\nu$ in place of $a(Y)\,dY$.

The defining equivariance relation gives, for
$\lambda\in\Wh_f(\pi)$, $u\in L_f$ and $v\in\pi$,
\[
 \lambda\!\left(\pi(u)v\right)=\theta_f(u)\lambda(v).
\]
With the Fourier convention
\eqref{eq:Fourier}: for $b\in C_c^\infty(\g)$,
\begin{equation}
\widehat{(\psi(B(f,\cdot))b)}(\xi)
 =\int_{\g}b(Y)\psi(-B(\xi-f,Y))\,dY
 =\widehat b(\xi-f).
\end{equation}
The same identity holds for finite compactly supported measures.
Thus multiplication by $\psi(B(f,Y))$ translates the support of
the Fourier transform by $f$.

For a lattice $M$ in a subspace of $\g$, its annihilator is
\begin{equation}
\label{eqn: lattice-perp}
 M^\perp=
 \{\xi\in\g:\psi(B(\xi,Y))=1\ \text{for every }Y\in M\}.
\end{equation}
When $M\subset\g_2$, the notation $M^\perp_2$ means its annihilator
inside $\g_{-2}$.
All annihilators are taken with respect to the fixed character $\psi$.

For a function $a$ supported in a sufficiently small neighborhood
of the identity, define
\[
 a^\natural(Y)=
 \begin{cases}
 a(\exp Y),&Y\text{ in the corresponding logarithmic neighborhood},\\
 0,&\text{otherwise}.
 \end{cases}
\]
For $a\in C_c^\infty(G)$ and a smooth representation $\pi$ of $G$,
define the integrated action by
\[
 \pi(a)=\int_G a(g)\pi(g)\,dg.
\]
More generally, for a finite compactly supported complex measure
$\nu$ on $G$, define
\[
 \pi(\nu)v=\int_G\pi(g)v\,d\nu(g)\qquad(v\in\pi).
\]
The orbit map $g\mapsto\pi(g)v$ is locally constant and takes only
finitely many values on the compact support of $\nu$, so this
integral is a finite sum.

On a compact group, by \emph{normalized Haar measure} we mean Haar
measure of total mass one. For a compact subgroup $H$, let $\mu_H$
denote its normalized Haar measure, viewed as a measure on $G$.
Let $\delta_g$ denote the point
mass at $g\in G$.

Define for $a, b \in C_c^{\infty}(G)$, 
\[
 (a*b)(g)=\int a(x)b(x^{-1}g)\,dx,
 \qquad a^*(g)=\overline{a(g^{-1})}.
\]
We further define convolution of measures by pushforward of the product measure by
group multiplication \cite[\S2.5, p.~54, equation~(2.34)]{Folland}.

A function $a$ is selfadjoint if $a^*=a$.
A function or measure is an idempotent if its convolution square
equals itself.

\subsection{The local character expansion and the known implication}
\label{subsec: local-character}

For each rational nilpotent orbit $\OO$, fix a positive invariant
orbital measure $\mu_{\OO}$. Existence and local finiteness of these
measures are supplied by \cite[Lemmas~43--44 and Theorem~45]{McNinch}.

For an admissible smooth representation $\pi$ of $G$ of finite length,
define its character distribution by
$\Theta_\pi(a)=\Tr\pi(a)$ for $a\in C_c^\infty(G)$.
Harish-Chandra's local character expansion
\cite[Theorem~16.2]{HC} gives a neighborhood
$\mathcal U$ of the identity on which
\begin{equation}\label{eq:character-expansion}
 \Theta_\pi(a)=
 \sum_{\OO}c_{\OO}(\pi)
        \int_{\OO}\widehat{a^\natural}(\xi)\,d\mu_{\OO}(\xi),
 \qquad a\in C_c^\infty(\mathcal U).
\end{equation}
To obtain the expansion for finite length, let
$\pi_1,\ldots,\pi_r$ be the irreducible factors of a composition
series of $\pi$, counted with multiplicity.
For a test function invariant on both sides under a compact open
subgroup $I$, its trace is computed on the finite-dimensional
space $\pi^I$. Taking $I$-invariants is exact by averaging, so
additivity of the trace gives additivity of the character.
Intersecting the finitely many neighborhoods for the expansions
of the $\pi_j$ therefore gives
\[
 \Theta_\pi=\sum_{j=1}^r\Theta_{\pi_j},
 \qquad
 c_{\OO}(\pi)=\sum_{j=1}^r c_{\OO}(\pi_j).
\]
Each measure $\mu_{\OO}$ has full support, so every nonempty open
subset of the orbit has positive measure. For every rational
nilpotent orbit $\OO$, \cite[Theorem~45]{McNinch} gives absolute
convergence:
\begin{equation}\label{eq:orbital-integrability}
 \int_{\OO}|b(\xi)|\,d\mu_{\OO}(\xi)<\infty
 \qquad\bigl(b\in C_c^\infty(\g)\bigr).
\end{equation}
There are finitely many rational nilpotent orbits for $G$
\cite[Theorem~40]{McNinch}. 
The coefficients $c_{\OO}(\pi)$ are taken with the Fourier convention
\eqref{eq:Fourier}.

In a give  a fixed geometric orbit, there are only finitely many rational
orbits in it, which are open and closed in its $F$-points
\cite[Theorem~40 and Lemma~44]{McNinch}. 

\begin{lemma}
The analytic closure relation defines a partial order on the finite
set of rational nilpotent orbits. 
\end{lemma}

\begin{proof}
Reflexivity and transitivity follow
from the properties of closure. 
Suppose
$\OO_1\subset\clan{\OO_2}$ with $\OO_1\ne\OO_2$.
Then $\OO_1$ and $\OO_2$ belong to distinct geometric orbits,
and $\OO_1$ lies in the Zariski boundary of the geometric orbit
containing $\OO_2$. Consequently $\dim\OO_1<\dim\OO_2$
\cite[Proposition~1.65 and the proof of Proposition~1.66]{MilneAG}.
This proves antisymmetry.
\end{proof}

One direction of our main result follows from \cite[Theorem~4.3.1(i)]{GGS17}, 
which attributes to  M{\oe}glin
and Waldspurger
\cite[Proposition~I.11, Theorem~I.16 and Corollary~I.17]{MW87}
and Varma \cite[Proposition~1 and Theorem~1]{Varma14} for residue characteristic two extension.

\begin{theorem}
\label{thm:published-forward}
Let $\pi$ be an admissible smooth complex representation of $G$
of finite length, and let $\OO$ be a rational nilpotent orbit. Then
\begin{equation}\label{eq:published-forward}
 \Wh_{\OO}(\pi)\ne0
 \quad\Longrightarrow\quad
 \OO\subset\operatorname{WF}(\pi).
\end{equation}
\end{theorem}

Only the converse of Theorem~\ref{thm:published-forward} remains
to be proved.

\section{Neutral vanishing and compactly supported measures}
\label{sec:neutral}

Fix a nonzero nilpotent $f$ and the neutral data from
Section~\ref{subsec:neutral-model}. Put
\[
 A=N_2/N_3\simeq(\g_2,+),\qquad
 G_0=Z_G(h)(F),\qquad \Omega_f=\Ad(G_0)f\subset\g_{-2}.
\]
The group $N_3$ is normal in $G_0N_1$, and $A$ is central in
$N_1/N_3$. For $u\in N_2$, write $\bar u=(\log u)_2$ for its
image in $A$. The smooth characters of $A$ are
\[
 \theta_\xi(a)=\psi(-B(\xi,a))
 \qquad(\xi\in\g_{-2},\ a\in A).
\]

\subsection{Smooth representations and twisted coinvariants}

\begin{lemma}\label{lem:additive-fibers}
Let $T$ be a locally compact abelian group that is a union of
compact subgroups, and let $W$ be a smooth representation of $T$.
Then $W=0$ if and only if
\[
 W_{T,\chi}=0
 \qquad\text{for every smooth character }\chi:T\to\C^\times.
\]
\end{lemma}

\begin{proof}
The main point is that the functor $X\mapsto X_{T,\chi}$ is exact,
by compact averaging and passage to the directed union of compact
subgroups (see \cite[Proposition~2.35(b)]{BZ76}). Consequently,
if $W_{T,\chi}=0$, then $X_{T,\chi}=0$ for every subquotient $X$ of $W$.

Note that every nonzero finitely generated submodule  has an irreducible
quotient.  Thus, if $W\ne0$, one of its subquotients is a
character $\chi$, contradicting $W_{T,\chi}=0$ for every $\chi$.
The converse is immediate.
\end{proof}

\begin{corollary}\label{cor:compact-invariants}
Let $T$ be a locally compact abelian group that is a union of
compact subgroups, and let $W$ be a smooth representation of $T$. Let $J\subset T$ be a compact open subgroup. Then $W^J=0$ if
and only if
\[
 W_{T,\chi}=0
 \qquad\text{for every smooth character }\chi:T\to\C^\times
       \text{ trivial on }J.
\]
\end{corollary}

\begin{proof}
Apply Lemma~\ref{lem:additive-fibers} to $W^J$, regarded as a
representation of $T/J$, which is again a locally compact group and is a union of compact subgroups.
\end{proof}

\begin{corollary}\label{cor:additive-compact-invariants}
Let $E$ be a finite-dimensional $F$-space, let $J\subset E$ be a
compact open subgroup, and let $W$ be a smooth representation of
the additive group $E$. Then $W^J=0$ if and only if
\[
 W_{E,\chi}=0
 \qquad\text{for every smooth character }\chi:E\to\C^\times
       \text{ trivial on }J.
\]
\end{corollary}

\begin{proof}
Corollary~\ref{cor:compact-invariants} applies because a
finite-dimensional vector space over $F$ is a union of compact
open additive subgroups.
\end{proof}

For a smooth representation $\pi$ of $G$, put $V=\pi_{N_3}$,
and write $[w]$ for the image of $w\in\pi$ in $V$.
The action of $A$ on $V$ is
\[
 a\cdot[w]=[\pi(\exp a)w]\qquad(a\in\g_2,\ w\in\pi),
\]
which gives 
$\pi_{N_2,\theta_f}\simeq V_{A,\theta_f}$.

\begin{lemma}\label{lem:neutral-reduction}
The following conditions are equivalent:
\begin{enumerate}
\renewcommand{\labelenumi}{(\theenumi)}
\item $\pi_{N_2,\theta_f}=0$;
\item $\pi_{L_f,\theta_f}=0$;
\item $\Wh_f(\pi)=0$.
\end{enumerate}
\end{lemma}

\begin{proof}
The first condition implies the second, and the second is
equivalent to the third by \eqref{eq:Wh-definition}.
It suffices to prove that (2) implies (1).
If $\g_1=0$, then $L_f=N_2$ and the assertion is immediate.

Let
\[
 \mathfrak n_2^f=\{X\in\g_2:B(f,X)=0\},\qquad
 U_f=\exp(\mathfrak n_2^f+\g_{\ge3}).
\]
The subspace $\mathfrak n_2^f$ has codimension one in $\g_2$.
The quotient
\[
 \g_{\ge1}/(\mathfrak n_2^f+\g_{\ge3})
\]
is therefore a Heisenberg Lie algebra with center
$\g_2/\mathfrak n_2^f$, and $\pi_{N_2,\theta_f}$ is a smooth
representation of $N_1/U_f$ with central character $\theta_f$.
Characters of $L_f/U_f$ with a different central restriction
have zero coinvariants on $\pi_{N_2,\theta_f}$.
Restricting to the abelian subgroup $L_f/U_f$ and applying
Lemma~\ref{lem:additive-fibers}, we obtain
\[
 \pi_{N_2,\theta_f}=0
 \quad\Longleftrightarrow\quad
 \pi_{L_f,\chi}=0
\]
for every character $\chi$ of $L_f/U_f$ extending
$\theta_f|_{N_2/U_f}$. By nondegeneracy of $\omega_f$,
conjugation by $N_1/U_f$ acts transitively on these characters.
The operators $\pi(n)$, for $n\in N_1$, identify the corresponding
coinvariant quotients. Thus (2) implies (1).
\end{proof}

\subsection{Lifting an indicator function to a unipotent subgroup}

Fix Haar measure $da$ on $A$. Choose Haar measure $d\xi$ on
$\g_{-2}$ so that, for $q\in C_c^\infty(A)$,
\begin{equation}\label{eq:weight-two-Fourier-inversion}
 \widehat q(\xi)=\int_A q(a)\psi(-B(\xi,a))\,da,\qquad
 q(a)=\int_{\g_{-2}}\widehat q(\xi)\psi(B(\xi,a))\,d\xi.
\end{equation}
Write $\vol_A$ and $\vol_{\g_{-2}}$ for the corresponding
additive volumes.

\begin{lemma}\label{lem:annihilator-volume}
For every lattice $R\subset A$, with annihilator $R^\perp_2$
defined in \eqref{eqn: lattice-perp}, one has
\[
 \vol_A(R)\vol_{\g_{-2}}(R^\perp_2)=1.
\]
\end{lemma}

\begin{proof}
Character orthogonality gives, for $\xi\in\g_{-2}$,
\[
 \widehat{\one_R}(\xi)
 =\int_R\psi(-B(\xi,a))\,da
 =\vol_A(R)\one_{R^\perp_2}(\xi).
\]
Fourier inversion at zero therefore gives
\[
 1=\one_R(0)
 =\int_{\g_{-2}}\widehat{\one_R}(\xi)\,d\xi
 =\vol_A(R)\vol_{\g_{-2}}(R^\perp_2).
\]
\end{proof}

For a compact open $C\subset\g_{-2}$, define
\begin{equation}\label{eq:compact-Fourier-action}
 q_C(a)=\int_C\psi(B(\xi,a))\,d\xi,\qquad
 \one_C\cdot v=\int_Aq_C(a)(a\cdot v)\,da\quad(v\in V).
\end{equation}
The function $q_C$ is compactly supported and locally constant.

\begin{lemma}\label{prop:open-spectrum}
The rational orbit $\Omega_f$ is open in $\g_{-2}$.
If $\Wh_f(\pi)=0$, then $\one_C\cdot V=0$ for every compact
open $C\subset\Omega_f$.
\end{lemma}

\begin{proof}
The orbit map $G_0\to\g_{-2}$, $g\mapsto\Ad(g)f$, has
surjective differential $X\mapsto[X,f]$ at the identity
\cite[Lemma~2.2.1(ii) and the paragraph following it]{GGS17}.
The inverse function theorem
\cite[Part~II, Chapter~III, \S9, Theorem~2]{SerreLALG}
therefore implies that $\Omega_f$ is open.

Now suppose $\Wh_f(\pi)=0$. The group $G_0$ acts on $V$, and
$\theta_{\Ad(g)f}(a)=\theta_f(\Ad(g^{-1})a)$.
Hence conjugation and Lemma~\ref{lem:neutral-reduction} give
\[
 V_{A,\theta_\xi}\simeq V_{A,\theta_f}
 \simeq\pi_{N_2,\theta_f}=0
 \qquad(\xi\in\Omega_f).
\]
Choose a lattice $R\subset A$ containing $\supp(q_C)$.
The annihilator $R^\perp_2$ is defined in
\eqref{eqn: lattice-perp}. Fourier inversion shows that
$\one_C$ is invariant under translation by $R^\perp_2$, so
\begin{equation}\label{eq:compact-coset-average}
 C=\coprod_{\nu=1}^r(\xi_\nu+R^\perp_2),\qquad
 q_C(a)=\frac{\one_R(a)}{\vol_A(R)}
              \sum_{\nu=1}^r\theta_{\xi_\nu}(a)^{-1}.
\end{equation}
The normalization follows from Lemma~\ref{lem:annihilator-volume}.
The characters of $A$ whose restrictions to $R$ equal
$\theta_{\xi_\nu}|_R$ are precisely $\theta_\eta$ with
$\eta\in\xi_\nu+R^\perp_2\subset C$.
Corollary~\ref{cor:additive-compact-invariants}, applied to
$V\otimes\theta_{\xi_\nu}^{-1}$, therefore gives
\[
 (V\otimes\theta_{\xi_\nu}^{-1})^R=0,\qquad
 \frac1{\vol_A(R)}
 \int_R\theta_{\xi_\nu}(a)^{-1}(a\cdot v)\,da=0.
\]
Summing over $\nu$ proves the assertion.
\end{proof}

Let $C\subset\Omega_f$ be compact open, and choose a lattice
$R\subset A$ containing $\supp(q_C)$.
A lattice $M\subset\g$ is \emph{graded} if
$M=\bigoplus_jM_j$, where $M_j=M\cap\g_j$; write
$M_{\ge i}=\bigoplus_{j\ge i}M_j$.
A compact open subgroup of $N_2$ is graded if its logarithm
is a graded lattice. Fix such a subgroup $U$ with image $R$ in $A$.
Choose Haar measures on $N_2$ and $N_3$ whose quotient measure
on $A$ is the fixed measure $da$.
For every graded compact open subgroup $K\subset N_2$ with
image $R$, define a measure supported on $K$ by
\begin{equation}\label{eq:subgroup-projector}
 d\phi_{C,K}(k)
 =\frac{q_C(\bar k)}{\vol(K\cap N_3)}\,dk
 \qquad(k\in K).
\end{equation}
Here $dk$ is Haar measure on $N_2$ restricted to $K$, and the
denominator uses Haar measure on $N_3$.

\begin{proposition}\label{lem:killing}
Suppose that $\Wh_f(\pi)=0$. For every finite-dimensional
subspace $D\subset\pi$, there is a graded compact open subgroup
$K\supset U$ of $N_2$ whose image in $A$ is exactly $R$ and
such that $\pi(\phi_{C,K})D=0$.
\end{proposition}

\begin{proof}
Put $S_R=\{u\in N_2:\bar u\in R\}$. We first show that every
compact subset $T\subset S_R$ is contained in a graded compact
open subgroup $K\supset U$ with image $R$ in $A$.
The Baker--Campbell--Hausdorff formula
\cite[14.35 and Proposition~14.36]{MilneAG} gives
\begin{equation}\label{eq:graded-BCH}
 \log(\exp X\exp Y)_j
 =X_j+Y_j+P_j(X_{<j},Y_{<j}),\qquad
 X_{<j}=\sum_{2\le i<j}X_i,
\end{equation}
for $X,Y\in\g_{\ge2}$, where $P_j$ contains the nonlinear
terms of weight $j$.
Set $\kappa_2=R$. For $j\ge3$, successively choose lattices
$\kappa_j\subset\g_j$ containing the weight-$j$ components of
$\log U$ and $\log T$, and $P_j(\kappa_{<j},\kappa_{<j})$,
where $\kappa_{<j}=\bigoplus_{2\le i<j}\kappa_i$.
These sets are compact and involve only previously chosen lattices.
Thus $\kappa=\bigoplus_{j\ge2}\kappa_j$ is closed under the
Baker--Campbell--Hausdorff product and negation.
Consequently $K=\exp\kappa$ is a graded compact open subgroup
containing $U\cup T$, with image $\kappa_2=R$ in $A$.
These graded compact open subgroups form a directed family,
and every compact subset of $S_R$ is contained in one of them.

Using the decomposition \eqref{eq:compact-coset-average}, define
$\chi_\nu(u)=\theta_{\xi_\nu}(\bar u)$ on $S_R$ for
$1\le\nu\le r$.
As in the proof of Lemma~\ref{prop:open-spectrum},
Corollary~\ref{cor:additive-compact-invariants} gives
\[
 \pi_{S_R,\chi_\nu}
 \simeq (V\otimes\theta_{\xi_\nu}^{-1})_R
 \simeq (V\otimes\theta_{\xi_\nu}^{-1})^R=0.
\]
We recall the argument of the Jacquet--Langlands lemma
\cite[Lemma~2.33]{BZ76}. Choose a basis $v_1,\ldots,v_d$
of $D$. By the definition, 
$\pi_{S_R,\chi_\nu}=0$ implies that, for each $i$ and $\nu$,
there is a finite expression
\[
 v_i=\sum_{j=1}^{m_{i\nu}}
 \bigl(\pi(g_{i\nu j})w_{i\nu j}
       -\chi_\nu(g_{i\nu j})w_{i\nu j}\bigr),
\]
where $g_{i\nu j}\in S_R$ and $w_{i\nu j}\in\pi$.
Set
\[
 T=\bigcup_{i=1}^d\bigcup_{\nu=1}^r
   \{g_{i\nu j}:1\le j\le m_{i\nu}\}.
\]
The preceding construction
therefore gives a graded compact open subgroup $K\supset U\cup T$
with image exactly $R$ in $A$.
For $g\in K$ and $w\in\pi$, the substitution $k'=kg$ and right
invariance of $\mu_K$ give
\[
 \begin{aligned}
 \pi(\chi_\nu^{-1}\mu_K)\pi(g)w
 &=\int_K\chi_\nu(k)^{-1}\pi(kg)w\,d\mu_K(k)\\
 &=\chi_\nu(g)\int_K\chi_\nu(k')^{-1}\pi(k')w\,d\mu_K(k')\\
 &=\chi_\nu(g)\pi(\chi_\nu^{-1}\mu_K)w.
 \end{aligned}
\]
Since every $g_{i\nu j}$ lies in $K$, applying this identity to
each summand in the expression for $v_i$ gives
\[
 \pi(\chi_\nu^{-1}\mu_K)D=0\qquad(1\le\nu\le r).
\]

Finally, the quotient Haar normalization gives
$\vol(K)=\vol(K\cap N_3)\vol_A(R)$, and hence
\[
 \phi_{C,K}=\sum_{\nu=1}^r\chi_\nu^{-1}\mu_K.
\]
Each summand annihilates $D$, so $\pi(\phi_{C,K})D=0$.
\end{proof}

\begin{remark}
By character orthogonality, the measure $\phi_{C,K}$ is a
convolution idempotent, but this is not needed for us.
\end{remark}

The next two sections prepare the
convolution argument in Section~\ref{sec:positive}, which produces
a function supported in that neighborhood whose action on $\pi$ is zero.

\section{Rational Slodowy slices and the adapted lattices}
\label{sec:geometry}

\subsection{Orbit closures and a rational Slodowy slice}

Recall that $\g^z=\ker(\ad z)$ is the Lie algebra centralizer of
$z\in\g$. For the triple $(e,h,f)$ fixed in
Section~\ref{subsec:neutral-model}, the Slodowy slice is
\[
 S_f=f+\g^e,\qquad
 \g^e=\bigoplus_{j\ge0}(\g^e\cap\g_j).
\]
The second equality follows from the highest-weight decomposition
recalled in Section~\ref{subsec:neutral-model}.

\begin{lemma}\label{lem:slice-incidence}
For every rational nilpotent orbit $\OO'$,
\[
 f\in\clan{\OO'}
 \quad\Longleftrightarrow\quad
 \OO'\cap S_f\ne\varnothing.
\]
If $s=f+\sum_{j\ge0}X_j$ belongs to this intersection, with
$X_j\in\g^e\cap\g_j$, then for $u\in F^\times$,
\begin{equation}\label{eq:Kazhdan}
 \kappa_u(s)=u^2\Ad(\gamma_h(u))s
           =f+\sum_{j\ge0}u^{j+2}X_j
\end{equation}
belongs to the rational orbit of $s$ for every $u\ne0$
and tends to $f$ as $u\to0$.
\end{lemma}

\begin{proof}
Choose an $F$-linear complement $\mathfrak m$ to $\g^f$ in $\g$
and a neighborhood $\mathcal V$ of zero in $\mathfrak m$ on which
the exponential is defined. Consider
\[
 \Phi:\mathcal V\times\g^e\longrightarrow\g,
 \qquad \Phi(Y,Z)=\Ad(\exp Y)(f+Z).
\]
The $\mathfrak{sl}_2$ decomposition
\cite[Chapter~VIII, \S1, no.~2, Propositions~1--2]{Bourbaki789} gives
\[
 \g=[\g,f]\oplus\g^e,\qquad
 d\Phi_{(0,0)}(Y,Z)=[Y,f]+Z.
\]
Since $Y\mapsto[Y,f]$ identifies $\mathfrak m$ with $[\g,f]$,
the derivative is an isomorphism. By the $F$-analytic inverse
function theorem
\cite[Part~II, Chapter~III, \S9, Theorem~2]{SerreLALG}, $\Phi$ maps
a neighborhood of $(0,0)$ isomorphically onto a neighborhood
$\mathcal W$ of $f$.
If $f\in\clan{\OO'}$, choose $x\in\OO'\cap\mathcal W$ and write
$x=\Phi(Y,Z)$. Then
\[
 f+Z=\Ad(\exp(-Y))x\in\OO'\cap S_f.
\]

For the converse, let $s=f+\sum_{j\ge0}X_j\in\OO'\cap S_f$,
with $X_j\in\g^e\cap\g_j$. The cocharacter $\gamma_h$ acts on
$\g_j$ by multiplication by $u^j$. Since $f\in\g_{-2}$,
\[
 \begin{aligned}
 u^2\Ad(\gamma_h(u))s
 =u^2\left(u^{-2}f+\sum_{j\ge0}u^jX_j\right)=f+\sum_{j\ge0}u^{j+2}X_j,
 \end{aligned}
\]
which proves \eqref{eq:Kazhdan}.
\end{proof}

\subsection{A uniform bound for unipotent conjugation}

Choose a finite-dimensional rational $F$-representation
\[
 \iota:G\hookrightarrow\operatorname{GL}_F(V)
\]
that is a closed immersion; such a representation exists by
\cite[Theorem~4.9 and Corollary~4.10]{MilneAG}.
Identify $G$ and $\g$ with their images under $\iota$ and $d\iota$.
The cocharacter $\gamma_h$ gives the weight decompositions
\[
 V=\bigoplus_a V_a,\qquad
 V_a=\{v\in V:\iota(\gamma_h(u))v=u^av\text{ for all }u\in F^\times\},
\]
\[
 \End_F(V)_j=\{M:M(V_a)\subset V_{a+j}\text{ for every }a\},
 \qquad E_+=\bigoplus_{j>0}\End_F(V)_j.
\]
The induced topology on $G$ from the matrix space is its $F$-analytic topology. Analytic exponential and logarithm
agree with the matrix exponential and logarithm on sufficiently small
neighborhoods
\cite[Part~II, Chapter~V, \S7, Corollary~2 and the following
remark]{SerreLALG}.
Fix a nonarchimedean norm on the finite-dimensional matrix space and
use its restriction as the norm on the Lie algebra.
For $s\in\g$, let
$\End_F(V)^s=\{M\in\End_F(V):Ms=sM\}$ be the centralizer of $s$ in $\End_F(V)$.

\begin{lemma}\label{lem:matrix-transporter}
If $s\in S_f$, then $\End_F(V)^s\cap E_+=0$.
There are a neighborhood $W$ of $s$ in the Lie algebra and a constant $c>0$
such that for $x,y\in W, k\in N_3$ with $\Ad(k)x=y$ we have 
\begin{equation}\label{eq:matrix-bound}
 \|k-1\|\le c\|y-x\|
 \end{equation}
\end{lemma}

\begin{proof}
Write $s=f+X$, where $X\in \g^e$ is a sum of components of nonnegative weight.
Write $M_j$ for the weight-$j$ component of a matrix $M$.
Suppose $M\in E_+$ is nonzero and $[s,M]=0$. Let
$j=\min\{i:M_i\ne0\}>0$. Since $X$ has only nonnegative weights,
\[
 0=[s,M]_{j-2}=[f,M_j].
\]
The action of $f$ on the positive-weight spaces of the
$\mathfrak{sl}_2$-module $\End_F(V)$ is injective
\cite[Chapter~VIII, \S1, no.~2, Proposition~2 and its
corollary]{Bourbaki789}. Hence $M_j=0$, a contradiction.

The linear map
\[
 \mathcal L_{x,y}:E_+\longrightarrow\End_F(V),
 \qquad M\longmapsto Mx-yM
\]
is therefore injective at $(x, y)=(s,s)$.
Put $d=\dim_F E_+$.  Fix a basis of $V$ and choose a projection
$P:\End_F(V)\to F^d$ onto $d$ coordinates such that
$P\mathcal L_{s,s}:E_+\to F^d$ is invertible.
On $F^d$, set
\[
 \|(a_1,\ldots,a_d)\|_\infty=\max_{1\le i\le d}|a_i|_F.
\]
Use the fixed matrix norm on $\End_F(V)$ and its restriction to
$E_+$.
After shrinking a compact neighborhood $W$ of $s$, the operators
$P\mathcal L_{x,y}$ remain invertible for $x,y\in W$, and
their inverses depend continuously on $(x,y)$.
Since $P$ is fixed and $W$ is compact, there is a constant $c>0$
such that
\[
 \|(P\mathcal L_{x,y})^{-1}P\|_{\mathrm{op}}
 :=\sup_{0\ne Z\in\End_F(V)}
   \frac{\|(P\mathcal L_{x,y})^{-1}PZ\|}{\|Z\|}
 \le c
 \qquad(x,y\in W).
\]
Consequently,
\[
 \begin{aligned}
 M&=(P\mathcal L_{x,y})^{-1}P\mathcal L_{x,y}M,\\
 \|M\|&\le\|(P\mathcal L_{x,y})^{-1}P\|_{\mathrm{op}}\, 
             \|\mathcal L_{x,y}M\|
       \le c\|\mathcal L_{x,y}M\|.
 \end{aligned}
\]

For $k\in N_3$, write $k=\exp Y$ with $Y\in\g_{\ge3}$.
The matrix $Y$ raises weights on $V$, so
\[
 k-1=\sum_{r\ge1}\frac{Y^r}{r!}\in E_+,
\]
where the sum is finite. If $\Ad(k)x=y$, then $kx=yk$ and
\[
 \mathcal L_{x,y}(k-1)
 =(k-1)x-y(k-1)=kx-yk+y-x=y-x.
\]
Substituting $M=k-1$ into the norm estimate proves
\eqref{eq:matrix-bound}.
\end{proof}

In particular, the estimate gives
\[
 \|k-1\|\le c\max\{\|x-s\|,\|y-s\|\}.
\]
Thus every $k\in N_3$ with $\Ad(k)x=y$ tends to the identity as
$x,y\to s$.

\subsection{An adapted lattice}

Write $|\cdot|_F$ for the absolute value on $F$.
For a full $\mathcal O_F$-lattice $M$ with an $\mathcal O_F$-basis
$v_1,\ldots,v_d$,
the lattice norm is
$\|\sum_i a_iv_i\|_M=\max_i|a_i|_F$.

For $s\in S_f$, a full graded $\mathcal O_F$-lattice
$\Lambda=\bigoplus_j(\Lambda\cap\g_j)$ in $\g$ is called
\emph{adapted to $s$}, relative to the fixed triple $(e,h,f)$,
if the exponential map is defined on $\Lambda$,
$[s,\Lambda]\subset\Lambda$, and the following conditions hold:
\begin{equation}\label{eq:powerful-bracket}
 [\Lambda,\Lambda]\subset p\Lambda\quad(p\ne2),
 \qquad
 [\Lambda,\Lambda]\subset8\Lambda\quad(p=2),
\end{equation}
\begin{equation}\label{eq:depth-isometry}
 \|\Ad(\exp Y)s-s\|_\Lambda=\|Y\|_\Lambda,\qquad
 \Ad(\exp Y)\Lambda=\Lambda
 \quad(Y\in\Lambda_{\ge3}).
\end{equation}
Here $\Lambda_{\ge3}=\Lambda\cap\g_{\ge3}$.

\begin{lemma}\label{lem:adapted-lattice}
Let $s\in S_f$ be rational and nilpotent.
After replacing $s$ by a point of the same rational orbit
in $S_f$, every neighborhood of zero in $\g$ contains a lattice
adapted to $s$.
\end{lemma}

\begin{proof}
Write $s_0$ for the original point $s$ in the statement.
In each irreducible summand of the adjoint $\mathfrak{sl}_2$-module,
choose a highest-weight vector $v_0$ of weight $d$ and let
$v_r=(\ad f)^rv_0$ for $0\le r\le d$.
These vectors form a basis
\cite[Chapter~VIII, \S1, no.~2, Propositions~1--2]{Bourbaki789}, with
\[
 [h,v_r]=(d-2r)v_r,\qquad
 [f,v_r]=v_{r+1}\ (r<d),\qquad [f,v_d]=0.
\]
Let $\Gamma$ be the $\mathcal O_F$-span of the bases of all
summands, including the trivial summands in the center.
Then $\Gamma$ is full and graded, and $[f,\Gamma]\subset\Gamma$.
If $v_r$ has positive weight, then $r<d$; its image under $\ad f$
is the distinct basis vector $v_{r+1}$. Therefore for $Y\in\g_{>0}$,
\[
 \|[f,Y]\|_\Gamma=\|Y\|_\Gamma.
\]

Since $\kappa_u(s_0)-f\to0$ in \eqref{eq:Kazhdan}, choose
$u\in F^\times$ small enough that
$[\kappa_u(s_0)-f,\Gamma]\subset\varpi\Gamma$.
Fix such a $u$ and put
\[
 s:=\kappa_u(s_0)\in(G\cdot s_0)\cap S_f.
\]
Keep this new point $s$ fixed for the rest of the proof.
Then $[s,\Gamma]\subset\Gamma$, and
\[
 \|[s-f,Y]\|_\Gamma
 \le|\varpi|_F\|Y\|_\Gamma
 <\|[f,Y]\|_\Gamma
 \qquad(0\ne Y\in\g_{>0}).
\]
The ultrametric inequality now gives
\[
 \|[s,Y]\|_\Gamma
 =\|[f,Y]+[s-f,Y]\|_\Gamma=\|Y\|_\Gamma
 \qquad(Y\in\g_{>0}).
\]

With $s$ fixed, take $\Lambda=\varpi^b\Gamma$ with $b$ sufficiently large.
Since $\|Y\|_\Lambda=|\varpi|_F^{-b}\|Y\|_\Gamma$, the preceding
isometry also holds for $\|\cdot\|_\Lambda$.
Choose an integer $a\ge0$ such that
$[\Gamma,\Gamma]\subset\varpi^{-a}\Gamma$. Then
\[
 [\Lambda,\Lambda]
 =\varpi^{2b}[\Gamma,\Gamma]
 \subset\varpi^{b-a}\Lambda.
\]
Taking $b$ large enough makes $\Lambda$ as small as desired,
ensures convergence of the matrix exponential on $\Lambda$,
and gives \eqref{eq:powerful-bracket}.

By \cite[Chapter~III, \S4, no.~4, Corollary~3(i)]{Bourbaki123}, we can choose 
$\Lambda$ so that for any $Y\in\Lambda$,
\begin{equation}
\label{eqn: ad-commute-exp}
\Ad(\exp Y)=\exp(\ad Y).
\end{equation}

For $Y\in\Lambda_{\ge3}$, the above identity gives
\begin{equation}
\label{eqn: adjoint-exponential}
 \Ad(\exp Y)s-s
 =[Y,s]+\sum_{r\ge1}\frac{(\ad Y)^r[Y,s]}{(r+1)!}.
\end{equation}
The first term has lattice norm $\|Y\|_\Lambda$.
For $r\ge1$, the bracket inclusions give
\[
 \left\|\frac{(\ad Y)^r[Y,s]}{(r+1)!}\right\|_\Lambda
 \le
 \begin{cases}
 \left|\dfrac{p^r}{(r+1)!}\right|_F\|Y\|_\Lambda,&p\ne2,\\[4pt]
 \left|\dfrac{8^r}{(r+1)!}\right|_F\|Y\|_\Lambda,&p=2.
 \end{cases}
\]
Let $v_p$ denote the valuation on $\mathbb Q$ with $v_p(p)=1$.
Counting the multiples of each power of $p$ in $(r+1)!$ gives
\begin{equation}\label{eq:factorial-valuation}
 v_p((r+1)!)=\sum_{a\ge1}\left\lfloor\frac{r+1}{p^a}\right\rfloor
 \le\frac{r}{p-1}.
\end{equation}
Thus both scalar factors in the preceding estimate have absolute
value less than one. The ultrametric inequality proves
$\|\Ad(\exp Y)s-s\|_\Lambda=\|Y\|_\Lambda$; the equality also holds
for $Y=0$.

For every $Z\in\Lambda$ and $r\ge1$, the same bracket estimates give
\[
 \frac{(\ad Z)^r}{r!}\Lambda
 \subset
 \begin{cases}
 \dfrac{p^r}{r!}\Lambda,&p\ne2,\\[4pt]
 \dfrac{8^r}{r!}\Lambda,&p=2
 \end{cases}
 \subset\Lambda.
\]
These scalar factors tend to zero as $r\to\infty$. Hence
\[
 \Ad(\exp(\pm Z))\Lambda
 =\exp(\pm\ad Z)\Lambda\subset\Lambda,
 \qquad \Ad(\exp Z)\Lambda=\Lambda.
\]
In particular, lattice preservation holds for $Z\in\Lambda_{\ge3}$,
completing \eqref{eq:depth-isometry}.
\end{proof}

For a full lattice $M\subset\g$, recall that its annihilator is
\[
 M^\perp=\{\xi\in\g:\psi(B(\xi,X))=1\text{ for all }X\in M\}.
\]
The invariance of $B$ gives
$(\Ad(q)M)^\perp=\Ad(q)(M^\perp)$ for $q\in G$.

\begin{proposition}\label{prop:exact-transporter}
Let $s\in S_f$ and $\Lambda$ be a lattice adapted to $s$.
For every sufficiently large fixed integer $m$, let $ \Lambda_m=\varpi^m\Lambda$ and
\[
Q=\exp\Lambda_m,
  \qquad L=\alpha\Lambda_m, 
   \qquad  \alpha=\begin{cases}1,&p\ne2,\\2,&p=2.\end{cases}
\]
For all sufficiently small  $t\in F^\times$, the compact open set
\begin{equation}\label{eq:tube}
 \mathcal T_t=\Ad(Q)(t^{-2}s)+L^\perp
\end{equation}
satisfies the following equality of subsets of $N_3$:
\begin{equation}\label{eq:exact-transporter}
 \{k\in N_3:\mathcal T_t\cap\Ad(k)\mathcal T_t
                         \ne\varnothing\}
 =Q\cap N_3.
\end{equation}
\end{proposition}

\begin{proof}
Let $W$ and $c$ be as in Lemma~\ref{lem:matrix-transporter}.
By continuity of the logarithm, choose $\varepsilon>0$ such that
\[
 k\in N_3,\quad\|k-1\|<\varepsilon
 \quad\Longrightarrow\quad \log k\in\Lambda_{\ge3}.
\]
Put $C_\Lambda=\max_{Z\in\Lambda}\|Z\|$ and choose $m$ so large that
\[
 s+\varpi^m\Lambda\subset W,\qquad
 c|\varpi|_F^m C_\Lambda<\varepsilon.
\]
Suppose $x=s+a$, $y=s+b$ with $a,b\in\varpi^m\Lambda$, and
$\Ad(k)x=y$ for $k\in N_3$. Then
\[
 \|k-1\|\le c\|y-x\|
 \le c|\varpi|_F^m C_\Lambda<\varepsilon.
\]
Let $k=\exp Y$ with $Y\in\Lambda_{\ge3}$. By
\eqref{eq:depth-isometry}, $\Ad(k)\Lambda=\Lambda$ and
\[
 \begin{aligned}
 \|Y\|_\Lambda
 =\|\Ad(k)s-s\|_\Lambda=\|b-\Ad(k)a\|_\Lambda\le|\varpi|_F^m.
 \end{aligned}
\]
Therefore $Y\in\varpi^m\Lambda_{\ge3}$ and $k\in Q\cap N_3$.

As in \eqref{eqn: ad-commute-exp}, we choose $m$ so that for any $Z\in\Lambda_m$ we have 
$\Ad(\exp Z)=\exp(\ad Z)$.
The bracket inclusions
\eqref{eq:powerful-bracket} and factorial estimate
\eqref{eq:factorial-valuation}, together with
$[s,\Lambda_m]\subset\Lambda_m$, give the following convergent
series in $\Lambda_m$:
\[
 \Ad(\exp Z)s-s
 =\sum_{r\ge0}\frac{(\ad Z)^r[Z,s]}{(r+1)!}\in\Lambda_m.
\]
Thus $\Ad(Q)s\subset s+\Lambda_m$.
Since $L^\perp$ is compact, take $t\in F^\times$ small enough that
\[
 |t|_F^2\max_{\xi\in L^\perp}\|\xi\|_\Lambda\le|\varpi|_F^m.
\]
It follows that
\[
 t^2\mathcal T_t
 =\Ad(Q)s+t^2L^\perp
 \subset s+\varpi^m\Lambda.
\]
If $x,y\in\mathcal T_t$ and $\Ad(k)x=y$ for $k\in N_3$, then
\[
 t^2x,t^2y\in s+\varpi^m\Lambda,\qquad
 \Ad(k)(t^2x)=t^2y.
\]
The preceding bound therefore gives $k\in Q\cap N_3$.

Conversely, the lattice-preservation calculation in
Lemma~\ref{lem:adapted-lattice} gives $\Ad(q)\Lambda_m=\Lambda_m$
for every $q\in Q$. Hence $\Ad(q)L=L$ and $\Ad(q)L^\perp=L^\perp$.
Consequently
\[
 \Ad(q)\mathcal T_t
 =\Ad(qQ)(t^{-2}s)+\Ad(q)L^\perp
 =\Ad(Q)(t^{-2}s)+L^\perp=\mathcal T_t.
\]
Since $\mathcal T_t$ is nonempty, every $q\in Q\cap N_3$ belongs
to the left-hand side of \eqref{eq:exact-transporter}.
\end{proof}

\section{Central idempotents and their convolution with central measures}
\label{sec:compact}

We construct central idempotents on small compact open subgroups and
study their convolution with central measures, which will be used in 
later constructions.

\subsection{Haar measure under the exponential maps}
\label{subsec:exponential-haar}

The analytic exponential map has derivative
equal to the identity at zero
\cite[Part~II, Chapter~V, \S7, Corollary~2]{SerreLALG}.
Choose Haar measure $dg$ on $G$ so that the Jacobian determinant
$J(Y)$ of the exponential map, relative to the chosen additive
Haar measure $dY$, satisfies $J(0)=1$.
On a sufficiently small neighborhood $\mathcal V$ of zero,
$|J(Y)|_F=1$. Thus, for a function $a$ supported in $\exp\mathcal V$,
\[
 \int_G a(g)\,dg
 =\int_{\mathcal V}a(\exp Y)|J(Y)|_F\,dY
 =\int_\g a^\natural(Y)\,dY.
\]

Choose a full $\mathcal O_F$-lattice $\Delta\subset\mathcal V$
satisfying
\[
 [\Delta,\Delta]\subset
 \begin{cases}
 p\Delta,&p\ne2,\\
 4\Delta,&p=2.
 \end{cases}
\]
Regarded as a finite free $\mathbb Z_p$-Lie algebra, $\Delta$ is
\emph{uniform} in the sense of \cite[Section~2.1]{BS}.
The Baker--Campbell--Hausdorff product identifies $\Delta$ with a
compact open subgroup $Q=\exp\Delta$ of $G$.
Moreover, $\Ad(q)\Delta=\Delta$ for $q\in Q$, and the subgroups
$\exp(p^r\Delta)$, $r\ge0$, are normal in $Q$ and form a
neighborhood basis of the identity
\cite[Theorem~9.8 and its proof]{DDMS}.
At $p=2$, the same result gives $Q^2=\exp(2\Delta)$, where $Q^2$
denotes the closed subgroup generated by the squares in $Q$.

For the lattice $\Delta$ fixed above, define
\[
 L=\alpha\Delta,\qquad
 \alpha=\begin{cases}1,&p\ne2,\\2,&p=2.\end{cases}
\]
For any $k\in G$ and $X,Y\in L$, compatibility of the matrix
exponential and logarithm with conjugation gives
\[
 \exp X=k\exp Yk^{-1}
 \quad\Longleftrightarrow\quad
 X=\log(k\exp Yk^{-1})=\Ad(k)Y.
\]
 Consequently
\begin{equation}\label{eq:intersection-log}
 \exp L\cap k\exp Lk^{-1}
 =\exp\bigl(L\cap\Ad(k)L\bigr).
\end{equation}
Put $M=L\cap\Ad(k)L$ and $J=\exp M$. Since $M\subset\Delta\subset
\mathcal V$, the chosen Haar normalization above gives
\[
 \int_J a(q)\,dq=\int_M a(\exp Y)\,dY,
 \qquad \vol(J)=\vol_\g(M),
\]
where $\vol_\g$ denotes volume with respect to $dY$.

\subsection{Central idempotents}

A function or measure on $Q$ is called \emph{central} if it is
invariant under conjugation by $Q$.
Take $Q=\exp\Delta$ as above and $L=\alpha\Delta$.
Thus $\exp L=Q$ at odd $p$ and $\exp L=Q^2$ at $p=2$.
Fix $z\in\g$ and define
\[
 \mathcal T=\Ad(Q)z+L^\perp,
\]
where $L^\perp$ is defined in \eqref{eqn: lattice-perp}.
Since $\Ad(Q)z$ is compact and $L^\perp$ is open, $\mathcal T$
is a finite disjoint union of cosets of $L^\perp$.
Define $e_{\mathcal T}$ on $Q$ by
\begin{equation}\label{eq:compact-projector-formula}
 e_{\mathcal T}(\exp Y)=\int_{\mathcal T}\psi(B(\xi,Y))\,d\xi
 \qquad(Y\in\Delta),
\end{equation}
and set $e_{\mathcal T}(g)=0$ for $g\notin Q$.

For odd $p$, the construction \eqref{eq:compact-projector-formula}
is the central idempotent obtained from Howe's character formula
\cite[Lemma~1.2 and Theorem~1.1]{Howe77}, under the hypotheses
stated there.

We first recall the convolution identity needed below.
With the group and additive Haar measures identified as in
Section~\ref{subsec:exponential-haar}, for central smooth functions $a,b$ on $Q$ and $p$ odd,  
\[
 (a*b)^\natural(X)
 =\int_\g a^\natural(Y)b^\natural(X-Y)\,dY.
\]
This follows from 
\cite[proof of Theorem~2.6]{BS}.

At $p=2$, the same identity
holds if both functions are supported on $Q^2$
\cite[Proposition~2.7(b)]{BS}, or if
$[\Delta,\Delta]\subset8\Delta$ and one function is supported
on $Q^2$ \cite[Proposition~2.7(a)]{BS}.

Taking Fourier transforms gives
\begin{equation}\label{eq:central-function-comparison}
 \widehat{(a*b)^\natural}
 =\widehat{a^\natural}\,\widehat{b^\natural}.
\end{equation}

\begin{proposition}\label{prop:compact-projector}
The function $e_{\mathcal T}\in C_c^\infty(G)$ is supported on $\exp L$,
is invariant under conjugation by $Q$, and satisfies
\[
 e_{\mathcal T}*e_{\mathcal T}=e_{\mathcal T},\qquad e_{\mathcal T}^*=e_{\mathcal T}.
\]
\end{proposition}

\begin{proof}
For $\eta,Y\in\g$, additive character orthogonality and the Haar
normalization for Fourier inversion give
\[
 \int_{\eta+L^\perp}\psi(B(\xi,Y))\,d\xi
 =\frac{\one_L(Y)}{\vol_\g(L)}\psi(B(\eta,Y)).
\]
Summing over the finitely many cosets in $\mathcal T$ shows that
the inverse Fourier transform of $\one_{\mathcal T}$ is supported
on $L$. Since $L\subset\Delta$, this inverse Fourier transform
is $e_{\mathcal T}^\natural$ on all of $\g$.
It follows from the definition that 
 $e_{\mathcal T}\in C_c^\infty(G)$ and
$\widehat{e_{\mathcal T}^\natural}=\one_{\mathcal T}$.

Since $Q$ preserves $L$, the set $\mathcal T$ is $\Ad(Q)$-invariant.
Invariance of $B$ gives $|\det\Ad(q)|_F=1$ for $q\in Q$.
Consequently, for $q\in Q$ and $Y\in\Delta$,
\[
 e_{\mathcal T}(q\exp(Y)q^{-1})
 =\int_{\mathcal T}\psi(B(\xi,\Ad(q)Y))\,d\xi
 =\int_{\mathcal T}\psi(B(\xi,Y))\,d\xi
 =e_{\mathcal T}(\exp Y).
\]
Thus $e_{\mathcal T}$ is central. It also follows from the definition
\[
 \overline{e_{\mathcal T}(\exp(-Y))}
 =\int_{\mathcal T}\psi(B(\xi,Y))\,d\xi
 =e_{\mathcal T}(\exp Y)\qquad(Y\in\Delta),
\]
so $e_{\mathcal T}$ is selfadjoint.

Apply \eqref{eq:central-function-comparison} with
$a=b=e_{\mathcal T}$ ( at $p=2$, both factors are supported on $Q^2$), 
we obtain
\[
 \widehat{(e_{\mathcal T}*e_{\mathcal T})^\natural}
 =\bigl(\widehat{e_{\mathcal T}^\natural}\bigr)^2
 =\one_{\mathcal T}^2
 =\one_{\mathcal T}
 =\widehat{e_{\mathcal T}^\natural}.
\]
Fourier inversion yields $e_{\mathcal T}*e_{\mathcal T}=e_{\mathcal T}$.
\end{proof}

Recall that $\delta_k$ denotes the point mass at $k\in G$.

\begin{lemma}\label{lem:compact-intertwining}
For the idempotent function $e_{\mathcal T}\in C_c^\infty(G)$, 
\[
 e_{\mathcal T}*\delta_k*e_{\mathcal T}\ne0
 \quad\Longrightarrow\quad
 \mathcal T\cap\Ad(k)\mathcal T\ne\varnothing.
\]
\end{lemma}

\begin{proof}
Since $\Ad(k)$ preserves $B$, one has $|\det\Ad(k)|_F=1$.
Let
\[
 e_{\mathcal T}^k(x)=e_{\mathcal T}(k^{-1}xk)
 =\bigl(\delta_k*e_{\mathcal T}*\delta_{k^{-1}}\bigr)(x),
 \qquad a=e_{\mathcal T}*e_{\mathcal T}^k.
\]
The function $e_{\mathcal T}^k$ is again a selfadjoint idempotent, and
\[
 e_{\mathcal T}*\delta_k*e_{\mathcal T}=a*\delta_k.
\]
Thus the hypothesis implies $a\ne0$.

For $b,c\in C_c^\infty(G)$, by invariance of Haar measure under
inversion we have 
\[
 (b*c)(1)=\int_G b(x)c(x^{-1})\,dx=(c*b)(1).
\]
Denote $\|a\|_2^2=\int_G|a(x)|^2\,dx$ and put
\[
 J=\exp L\cap k\exp Lk^{-1},\qquad
 M=L\cap\Ad(k)L.
\]
Using the fact that $e_{\mathcal T}$ and $e_{\mathcal T}^k$ are idempotents, we obtain 
\[
 \begin{aligned}
 0<\|a\|_2^2
 =(a^**a)(1)
   =(e_{\mathcal T}^k*e_{\mathcal T}*e_{\mathcal T}^k)(1)=(e_{\mathcal T}*e_{\mathcal T}^k)(1)
   =\int_J e_{\mathcal T}(q)\overline{e_{\mathcal T}(k^{-1}qk)}\,dq.
 \end{aligned}
\]
The integral is over $J$ because $e_{\mathcal T}$ is supported on $\exp L$.
By \eqref{eq:intersection-log}, $J=\exp M$.
By
\eqref{eq:compact-projector-formula} and Fourier inversion, we obtain
\[
 \begin{aligned}
 \|a\|_2^2
 =\int_\g e_{\mathcal T}^\natural(Y)
   \overline{e_{\mathcal T}^\natural(\Ad(k^{-1})Y)}\,dY=\int_{\mathcal T}
   \overline{\widehat{e_{\mathcal T}^\natural}(\Ad(k^{-1})\xi)}\,d\xi=\int_{\mathcal T}\one_{\mathcal T}(\Ad(k^{-1})\xi)\,d\xi>0.
 \end{aligned}
\]
The positive integral implies
$\mathcal T\cap\Ad(k)\mathcal T\ne\varnothing$.
\end{proof}

\subsection{Convolution with a central measure}

We extend \eqref{eq:central-function-comparison} to convolution
with a central measure.

For a finite complex measure $\Phi$ on $Q$, let $\log_*\Phi$
denote its pushforward to $\Delta$ by the logarithm. Thus
\[
 \int_\Delta b(Y)\,d(\log_*\Phi)(Y)
 =\int_Q b(\log q)\,d\Phi(q)
 \qquad(b\in C^\infty(\Delta)).
\]

\begin{lemma}\label{lem:measure-density}
Let $\Phi$ be a finite complex measure on $Q$, and let
$I_r=\exp(p^r\Delta)$ for $r\ge0$. Write $dq$ for the restriction
of $dg$ to $Q$ and $\vol(I_r)=\int_{I_r}dq$. Define
\[
 \Phi_r(q)=\frac{1}{\vol(I_r)}
       \int_Q\one_{I_r}(x^{-1}q)\,d\Phi(x),\qquad q\in Q.
\]
Then $\Phi_r$ is smooth, and
\[
 d(\Phi*\mu_{I_r})(q)=\Phi_r(q)\,dq.
\]
For every smooth function $e$ on $Q$,
\[
 e*\Phi_r=(e*\Phi)*\mu_{I_r}.
\]
\end{lemma}

\begin{proof}
Recall that $d\mu_{I_r}(h)=dh/\vol(I_r)$ on $I_r$.
For every continuous function $v$ on $Q$, Haar invariance and
Fubini's theorem give
\[
 \begin{aligned}
 \int_Q v(q)\Phi_r(q)\,dq
 &=\frac{1}{\vol(I_r)}\int_Q\int_{xI_r}v(q)\,dq\,d\Phi(x)\\
 &=\int_Q\int_{I_r}v(xh)\,d\mu_{I_r}(h)\,d\Phi(x)\\
 &=\int_Q v(q)\,d(\Phi*\mu_{I_r})(q).
 \end{aligned}
\]
This proves the density formula. Since
$\Phi_r(qh)=\Phi_r(q)$ for $h\in I_r$, the density is smooth.
For $g\in Q$, the density formula gives
\[
 \begin{aligned}
 (e*\Phi_r)(g)
 &=\int_Q e(gq^{-1})\Phi_r(q)\,dq\\
 &=\int_Q\int_{I_r}e(gh^{-1}x^{-1})\,
       d\mu_{I_r}(h)\,d\Phi(x)\\
 &=\int_{I_r}(e*\Phi)(gh^{-1})\,d\mu_{I_r}(h)
 =((e*\Phi)*\mu_{I_r})(g).
 \end{aligned}
\]
\end{proof}

\begin{lemma}\label{lem:central-measure}
Let $e$ be a central smooth function on $Q$, and let $\Phi$
be a finite central complex measure on $Q$.
For odd $p$, or for $p=2$ with
$[\Delta,\Delta]\subset8\Delta$ and $\supp(e)\subset Q^2$,
the convolution $e*\Phi$ is a central smooth function and
\begin{equation}\label{eq:central-measure-comparison}
 \widehat{(e*\Phi)^\natural}
 =\widehat{e^\natural}\,\widehat{\log_*\Phi}.
\end{equation}
\end{lemma}

\begin{proof}
Let $I_r$ and $\Phi_r$ be as in Lemma~\ref{lem:measure-density}.
Normality of $I_r$ and conjugation invariance of $\Phi$ give
\[
 \Phi_r(uqu^{-1})=\Phi_r(q)\quad(u\in Q).
\]
Thus the smooth function $\Phi_r$ is central.
For sufficiently large $r$, the function $e$ is right
$I_r$-invariant, so
\[
 (e*\mu_{I_r})(q)
 =\frac{1}{\vol(I_r)}\int_{I_r}e(qh^{-1})\,dh=e(q).
\]
Lemma~\ref{lem:measure-density}, centrality of $\mu_{I_r}$,
and associativity give
\[
 e*\Phi_r=e*\Phi*\mu_{I_r}
 =(e*\mu_{I_r})*\Phi=e*\Phi.
\]
In particular, $e*\Phi$ is a central smooth function.

Fix $\xi\in\g$. Local constancy on compact $Q$ gives, for all
sufficiently large $r$,
\[
 \psi(-B(\xi,\log(qh)))=\psi(-B(\xi,\log q))
 \qquad(q\in Q,\ h\in I_r).
\]
Using the  compatibility of Haar measures with the exponential map,
\[
 \begin{aligned}
 \widehat{\Phi_r^\natural}(\xi)
 &=\int_Q\int_{I_r}
   \psi(-B(\xi,\log(qh)))\,d\mu_{I_r}(h)\,d\Phi(q)\\
 &=\int_Q\psi(-B(\xi,\log q))\,d\Phi(q)
 =\widehat{\log_*\Phi}(\xi).
 \end{aligned}
\]
Choose $r$ large enough for both $e*\Phi_r=e*\Phi$ and the above
Fourier identity to hold. Applying
\eqref{eq:central-function-comparison} to $e$ and $\Phi_r$ gives
\[
 \begin{aligned}
 \widehat{(e*\Phi)^\natural}(\xi)
 =\widehat{(e*\Phi_r)^\natural}(\xi)=\widehat{e^\natural}(\xi)\widehat{\Phi_r^\natural}(\xi)
 =\widehat{e^\natural}(\xi)\widehat{\log_*\Phi}(\xi).
 \end{aligned}
\]
This proves \eqref{eq:central-measure-comparison} for every $\xi$.
At $p=2$, the required support condition holds for the fixed
factor $e$.
\end{proof}

\section{Construction of a function with nonnegative Fourier transform}
\label{sec:positive}

Let $\pi$ be the admissible smooth representation of finite length
under consideration in Theorem~\ref{thm:main}.
Fix a nonzero nilpotent $f$ and a rational nilpotent $s\in S_f$.
After replacing the chosen point within its rational orbit if necessary,
choose $s,\Lambda$ as in Lemma~\ref{lem:adapted-lattice}.
Choose a sufficiently large fixed $m$ for
Proposition~\ref{prop:exact-transporter} and let
\[
 \Lambda_m=\varpi^m\Lambda,\qquad Q=\exp\Lambda_m,\qquad
 L=\alpha\Lambda_m,\qquad
 \alpha=\begin{cases}1,&p\ne2,\\2,&p=2.\end{cases}
\]
We also require $m$ to be large enough that $Q\subset\mathcal U$,
where \eqref{eq:character-expansion} holds.
For sufficiently small $t\in F^\times$, recall that
\[
 \mathcal T_t=\Ad(Q)(t^{-2}s)+L^\perp.
\]
Apply Proposition~\ref{prop:compact-projector} with $z=t^{-2}s$ and
$\Delta=\Lambda_m$. The resulting idempotent $e_{\mathcal T_t}$ satisfies
\[
\widehat{e_{\mathcal T_t}^\natural}=\one_{\mathcal T_t}.
\]

\subsection{Restricting the convolution to compact subgroups}

Recall that $A=N_2/N_3\simeq\g_2$,
$G_0=Z_G(h)(F)$, and $\Omega_f=G_0f$ is open in $\g_{-2}$.
Let
\[
 R=L\cap\g_2,\qquad
 C_t=t^{-2}f+R^\perp_2,\qquad U=\exp L\cap N_2.
\]
Here, in the notation of \eqref{eqn: lattice-perp},
\[
 R^\perp_2
 =\{\xi\in\g_{-2}:\psi(B(\xi,a))=1
       \text{ for every }a\in R\}
\]
is the annihilator of $R\subset\g_2$ in $\g_{-2}$.

\begin{lemma}\label{lem:projection-cutoff}
Let $p_{-2}:\g\to\g_{-2}$ be the projection onto the weight-$-2$
subspace for $h$. Then
\[
 p_{-2}(t^{-2}s+L^\perp)=C_t,
 \qquad
 \mathcal T_t=\bigcup_{q\in Q}\Ad(q)(t^{-2}s+L^\perp).
\]
Consequently, for every $\xi\in\mathcal T_t$, there is $q\in Q$ such that
\[
 p_{-2}\bigl(\Ad(q^{-1})\xi\bigr)\in C_t.
\]
\end{lemma}

\begin{proof}
For $x\in L^\perp$ and $a\in R$, one has
\[
 \psi(B(p_{-2}x,a))=\psi(B(x,a))=1,
\]
so $p_{-2}x\in R_2^\perp$. Conversely, if $\eta\in R_2^\perp$
and $Y=\sum_jY_j\in L$ is the weight decomposition of $Y$, then
$Y_2\in R$ and
\[
 \psi(B(\eta,Y))=\psi(B(\eta,Y_2))=1.
\]
Thus $\eta\in L^\perp$, and $p_{-2}(L^\perp)=R_2^\perp$.
Since $s-f\in\g^e\subset\g_{\ge0}$, one has $p_{-2}(s)=f$.
These identities give $p_{-2}(t^{-2}s+L^\perp)=C_t$.

For $q\in Q$, the equality $\Ad(q)L=L$ and invariance of $B$
give $\Ad(q)L^\perp=L^\perp$. Hence
\[
 \mathcal T_t
 =\bigcup_{q\in Q}\bigl(\Ad(q)(t^{-2}s)+L^\perp\bigr)
 =\bigcup_{q\in Q}\Ad(q)(t^{-2}s+L^\perp).
\]
For $\xi\in\mathcal T_t$, choose $q\in Q$ with
$\Ad(q^{-1})\xi\in t^{-2}s+L^\perp$ and apply the previous identity.
\end{proof}

Since $R^\perp_2$ is compact and $\Omega_f$ is open,
$f+t^2R^\perp_2\subset\Omega_f$ for sufficiently small $t$.
The cocharacter $\gamma_h$ takes values in $G_0$, and
\[
 \Ad(\gamma_h(t))(f+t^2R^\perp_2)
 =t^{-2}f+R^\perp_2=C_t\subset\Omega_f.
\]
For $a\in A$, character orthogonality gives
\[
 \begin{aligned}
 q_{C_t}(a)
 &=\psi(B(t^{-2}f,a))
   \int_{R^\perp_2}\psi(B(\eta,a))\,d\eta\\
 &=\psi(B(t^{-2}f,a))\vol_{\g_{-2}}(R^\perp_2)\one_R(a).
 \end{aligned}
\]
Using Lemma~\ref{lem:annihilator-volume}, we obtain
\begin{equation}\label{eq:lower-cutoff}
 q_{C_t}(a)=
 \frac{\one_R(a)}{\vol_A(R)}\psi(B(t^{-2}f,a)).
\end{equation}

\begin{lemma}\label{lem:localized-sandwich}
Assume $\Wh_f(\pi)=0$.
There is a graded compact open subgroup $K$ of $N_2$ containing $U$,
whose image under the quotient map
\[
 N_2\longrightarrow N_2/N_3\simeq\g_2,\qquad
 u\longmapsto(\log u)_2,
\]
is $R=L\cap \mathfrak g_2$, and which satisfies $\pi(\phi_{C_t,K}*e_{\mathcal T_t})=0$.
Here $(\log u)_2$ denotes the $\g_2$-component of $\log u$.
For $H=K\cap N_3$, $H'=H\cap Q$, and $K'=K\cap Q$,
\begin{equation}\label{eq:localized-sandwich}
 e_{\mathcal T_t}*\phi_{C_t,K}*e_{\mathcal T_t}
 =
 \frac{\vol(H')}{\vol(H)}
       e_{\mathcal T_t}*\phi_{C_t,K'}*e_{\mathcal T_t}.
\end{equation}
Thus the convolution on the right annihilates $\pi$ and is
supported on $Q$.
\end{lemma}

\begin{proof}
Choose a compact open subgroup $I\subset Q$ under which $e_{\mathcal T_t}$
is left invariant. For $i\in I$,
\[
 \pi(i)\pi(e_{\mathcal T_t})=\pi(\delta_i*e_{\mathcal T_t})=\pi(e_{\mathcal T_t}).
\]
By admissibility, 
\[
 D:=\pi(e_{\mathcal T_t})\pi\subset\pi^I,\qquad
 \dim D\le\dim\pi^I<\infty.
\]
Proposition~\ref{lem:killing} supplies a graded $K\supset U$ with
image $R$ in $A$ and
\[
 \pi(\phi_{C_t,K})D=0,\qquad
 \pi(\phi_{C_t,K}*e_{\mathcal T_t})=0.
\]

The lattice $\log K'=(\log K)\cap\Lambda_m$ is graded.
For every $a\in R$, it follows from definition that $\exp(a)\in U\subset K'$.
Denote by $\bar k\in R$  the image of $k\in K$ in $A$, we obtain
\[
 k\exp(-\bar k)\in H,\qquad
 k\in K'\ \Longrightarrow\ k\exp(-\bar k)\in H'.
\]
Thus $K=HU$, $K'=H'U$, and both groups have image $R$ in $A$.

For $h\in H$ and $a\in R$, centrality of $e_{\mathcal T_t}$ on $Q$ gives
\[
 e_{\mathcal T_t}*\delta_{h\exp(a)}*e_{\mathcal T_t}
 =(e_{\mathcal T_t}*\delta_h*e_{\mathcal T_t})*\delta_{\exp(a)}.
\]
By Lemma~\ref{lem:compact-intertwining} and
Proposition~\ref{prop:exact-transporter},
 \begin{equation}
 \label{eqn: vanishing-condition-K'}
 e_{\mathcal T_t}*\delta_h*e_{\mathcal T_t}\ne0
  \Longrightarrow\
 \mathcal T_t\cap\Ad(h)\mathcal T_t\ne\varnothing \Longrightarrow\ h\in Q\cap N_3.
 \end{equation}

Since $h\in H$, the last condition is $h\in H'$.

For each $a\in R$, the fiber of $K\to R$ over $a$ is
$H\exp(a)$, and the fiber of $K'\to R$ over $a$ is
$H'\exp(a)$. With the quotient Haar measures,
fixed in \eqref{eq:subgroup-projector}
\[
 \begin{aligned}
 e_{\mathcal T_t}*\phi_{C_t,K}*e_{\mathcal T_t}
 &=\frac{1}{\vol(H)}
   \int_R q_{C_t}(a)
   \int_H e_{\mathcal T_t}*\delta_{h\exp(a)}*e_{\mathcal T_t}\,dh\,da\\
 &=\frac{1}{\vol(H)}
   \int_R q_{C_t}(a)
   \int_{H'} e_{\mathcal T_t}*\delta_{h\exp(a)}*e_{\mathcal T_t}\,dh\,da\\
 &=\frac{\vol(H')}{\vol(H)}
   e_{\mathcal T_t}*\phi_{C_t,K'}*e_{\mathcal T_t}.
 \end{aligned}
\]
The second equality uses the vanishing \eqref{eqn: vanishing-condition-K'}  for $h\notin H'$.
Since  $ \frac{\vol(H')}{\vol(H)}=[H:H']^{-1}>0 $, 
applying $\pi$ gives
\[
 0=\pi(e_{\mathcal T_t})\pi(\phi_{C_t,K})\pi(e_{\mathcal T_t})
 =\frac{\vol(H')}{\vol(H)}
   \pi(e_{\mathcal T_t}*\phi_{C_t,K'}*e_{\mathcal T_t}).
\]
The last convolution therefore annihilates $\pi$, and its support
is contained in $Q\,K'\,Q=Q$.
\end{proof}

\subsection{Positivity of the Fourier transform}

For the subgroup $K'$ supplied by
Lemma~\ref{lem:localized-sandwich}, let $\phi'=\phi_{C_t,K'}$.
By \eqref{eq:subgroup-projector}, it is supported on
$K'$ and is given by
\[
 d\phi'(k)=\frac{q_{C_t}(\bar k)}{\vol(H')}\,dk
 \qquad(k\in K'),
\]
where $\bar k$ is the image of $k$ in $N_2/N_3\simeq\g_2$,
and $dk$ is the restriction to $K'$ of the chosen Haar measure
on $N_2$.

For $q\in Q$, let $\Ad(q)_*$ denote pushforward of measures under
the conjugation map $k\mapsto qkq^{-1}$.
Define
\[
 \Phi:=\int_Q\Ad(q)_*\phi'\,dq,
\]
where $dq$ is normalized Haar measure on $Q$.
This defines a central compactly supported
complex measure on $Q$.

The proof below obtains a nonnegative Fourier transform by averaging
an indicator function under conjugation. A related finite-field
formula is recalled by Barbasch and Moy in
\cite[Proposition~2.6]{BM97}: Lusztig expresses the Fourier
transform of a generalized Gelfand--Graev character as a positive
multiple of the number of group elements carrying a given point
into a Slodowy slice.

\begin{lemma}\label{lem:measure-positivity}
The Fourier transform of the pushforward of $\Phi$ under the
logarithm is nonnegative everywhere and strictly positive near $t^{-2}s$.
\end{lemma}

\begin{proof}
Let
 \[
\kappa'=\log K',\qquad
\kappa'_{\ge3}=\kappa'\cap\g_{\ge3}.
\]
 The lattice $\kappa'$ is graded, with
$\kappa'\cap\g_2=R$.

To compute $\widehat{\log_*\phi'}$, we first rewrite 
integration on $K'$ as integration on $\kappa'$.
The triangular formula \eqref{eq:graded-BCH} shows that, for fixed $X\in\g_{\ge i}$,
the map $Y\mapsto\log(\exp X\exp Y)$ has Jacobian determinant
one on $\g_{\ge i}$, for $i=2,3$. Thus the exponential map
carries additive Haar measure on $\g_{\ge i}$ to Haar measure on $N_i$.

Choose additive Haar measure $dZ$ on $\g_{\ge3}$ to correspond
under the exponential map to the chosen Haar measure on $N_3$.
We use $\vol_{\g_{\ge3}}(D)$ to denote the volume of $D\subseteq \g_{\ge3}$ with respect to $dZ$.
Since $\kappa'=R\oplus\kappa'_{\ge3}$ and
$H'=\exp(\kappa'_{\ge3})$, the normalization gives
\[
 \vol_{\g_{\ge3}}(\kappa'_{\ge3})=\vol(H').
\]
With the quotient Haar normalization in
\eqref{eq:subgroup-projector}, we therefore have, for
$v\in C^\infty(K')$,
\[
 \int_{K'}v(k)\,dk
 =\int_R\int_{\kappa'_{\ge3}}v(\exp(a+Z))\,dZ\,da.
\]
For $a\in R$ and $Z\in\kappa'_{\ge3}$, the image of
$k=\exp(a+Z)$ in $A\simeq\g_2$ is $\bar k=a$.
By \eqref{eq:lower-cutoff},
\[
 q_{C_t}(\bar k)=q_{C_t}(a)
 =\frac{1}{\vol_A(R)}\psi(B(t^{-2}f,a)).
\]
It follows that 
\[
 \begin{aligned}
 \widehat{\log_*\phi'}(\xi)
 &=\int_{K'}\psi(-B(\xi,\log k))\,d\phi'(k)\\
 &=\frac{1}{\vol(H')}
   \int_{K'}q_{C_t}(\bar k)\psi(-B(\xi,\log k))\,dk\\
 &=\frac{1}{\vol(H')}
   \int_R\int_{\kappa'_{\ge3}}
   q_{C_t}(a)\psi(-B(\xi,a+Z))\,dZ\,da\\
 &=\frac{1}{\vol_A(R)\vol(H')}
   \int_R\int_{\kappa'_{\ge3}}
   \psi\bigl(B(t^{-2}f,a)-B(\xi,a+Z)\bigr)\,dZ\,da.
 \end{aligned}
\]
Denote by $\xi_{-2}$  the component of $\xi$ in $\g_{-2}$.
Since $a\in\g_2$ and $B$ pairs only opposite weight spaces,
\[
 B(t^{-2}f,a)-B(\xi,a+Z)
 =B(t^{-2}f-\xi_{-2},a)-B(\xi,Z).
\]
By Fubini's theorem, we obtain
\[
 \begin{aligned}
 \widehat{\log_*\phi'}(\xi)
 &=\left(\frac{1}{\vol_A(R)}
       \int_R\psi(B(t^{-2}f-\xi_{-2},a))\,da\right)
   \left(\frac{1}{\vol(H')}
       \int_{\kappa'_{\ge3}}\psi(-B(\xi,Z))\,dZ\right).
 \end{aligned}
\]
Since $\vol_{\g_{\ge3}}(\kappa'_{\ge3})=\vol(H')$, both
integrals are normalized by the volume of their domain.
Character orthogonality makes each normalized integral one when
the corresponding character is trivial and zero otherwise. Hence
\begin{equation}\label{eq:subgroup-Fourier}
 \widehat{\log_*\phi'}(\xi)=
 \one_{C_t}(\xi_{-2})\,
 \one_{(\kappa'_{\ge3})^\perp}(\xi)\ge0.
\end{equation}
Since $s-f\in\g^e$ has only nonnegative weights,
\[
 (t^{-2}s)_{-2}=t^{-2}f,\qquad
 B(t^{-2}s,\kappa'_{\ge3})=0,\qquad
 \widehat{\log_*\phi'}(t^{-2}s)=1.
\]

Invariance of $B$ under conjugation and compatibility of logarithm with conjugation give
\[
 \begin{aligned}
 \widehat{\log_*\Phi}(\xi)
 &=\int_Q\int_{K'}
   \psi\bigl(-B(\xi,\Ad(q)\log k)\bigr)\,d\phi'(k)\,dq\\
 &=\int_Q\widehat{\log_*\phi'}(\Ad(q^{-1})\xi)\,dq\ge0.
 \end{aligned}
\]
Since the integrand at $\xi=t^{-2}s$ is locally constant in $q$, we can choose an open subgroup $Q_0\subset Q$
on which it is one. By normalization,  we have $\int_Qdq=1$, hence
\begin{equation}\label{eq:averaged-positive}
 \widehat{\log_*\Phi}(t^{-2}s)
 \ge\int_{Q_0}dq=[Q:Q_0]^{-1}>0.
\end{equation}

Finally, we have $\supp(\log_*\Phi)\subset\Lambda_m$. For
$\zeta\in\Lambda_m^\perp$,
\[
 \begin{aligned}
 \widehat{\log_*\Phi}(\xi+\zeta)
 =\int_{\Lambda_m}\psi(-B(\xi,Y))\psi(-B(\zeta,Y))
                       \,d(\log_*\Phi)(Y)=\widehat{\log_*\Phi}(\xi).
 \end{aligned}
\]
Thus $\widehat{\log_*\Phi}$ is strictly positive throughout
the open neighborhood $t^{-2}s+\Lambda_m^\perp$.
\end{proof}

\begin{proposition}\label{prop:positive-test}
If $\Wh_f(\pi)=0$, then for every sufficiently small
nonzero $t$ there is $a_t\in C_c^\infty(Q)$ such that
\[
 \pi(a_t)=0,\qquad
 \widehat{a_t^\natural}\ge0,\qquad
 \supp(\widehat{a_t^\natural})\subset\mathcal T_t.
\]
Moreover,  $\widehat{a_t^\natural}$ is strictly positive on an
open neighborhood of $t^{-2}s$.
\end{proposition}

\begin{proof}
Let $b=e_{\mathcal T_t}*\phi'*e_{\mathcal T_t}$.
Lemma~\ref{lem:localized-sandwich} gives $\pi(b)=0$.
Centrality of $e_{\mathcal T_t}$ gives
\[
 \begin{aligned}
 a_t:=\int_Q\Ad(q)_*b\,dq
     =e_{\mathcal T_t}*\Phi*e_{\mathcal T_t}=(e_{\mathcal T_t}*e_{\mathcal T_t})*\Phi
      =e_{\mathcal T_t}*\Phi.
 \end{aligned}
\]
Note that
\[
 \pi(a_t)
 =\int_Q\pi(q)\pi(b)\pi(q^{-1})\,dq=0.
\]

At $p=2$, the hypotheses of Lemma~\ref{lem:central-measure} hold:
\[
 [\Lambda_m,\Lambda_m]\subset8\Lambda_m,\qquad
 \supp(e_{\mathcal T_t})\subset Q^2.
\]
It follows that
$a_t\in C_c^\infty(Q)$ and
\begin{equation}\label{eq:positive-test-Fourier}
 \widehat{a_t^\natural}
 =\one_{\mathcal T_t}\widehat{\log_*\Phi}.
\end{equation}
Lemma~\ref{lem:measure-positivity} shows the right-hand side
nonnegative. Moreover,
\[
 L\subset\Lambda_m
 \ \Longrightarrow\ \Lambda_m^\perp\subset L^\perp,
 \qquad
 t^{-2}s+\Lambda_m^\perp\subset t^{-2}s+L^\perp\subset\mathcal T_t.
\]
For $\xi\in t^{-2}s+\Lambda_m^\perp$, the indicator in
\eqref{eq:positive-test-Fourier} is one.
The invariance under translation by $\Lambda_m^\perp$ proved in
Lemma~\ref{lem:measure-positivity} and the bound
\eqref{eq:averaged-positive} therefore give
\[
 \widehat{a_t^\natural}(\xi)
 =\widehat{\log_*\Phi}(\xi)
 =\widehat{\log_*\Phi}(t^{-2}s)>0.
\]
Equation~\eqref{eq:positive-test-Fourier} also gives the stated
containment of the support.
\end{proof}

The construction uses the same compact subgroup $Q$ for every
sufficiently small $t\in F^\times$.

\section{Isolating an orbital integral and proof ot the main theorem}
\label{sec:conclusion}

Recall that by definition a rational nilpotent orbit $\OO$ occurs in
$\operatorname{WFC}(\pi)$ precisely when the coefficient
$c_{\OO}(\pi)$  is nonzero. Recall that the wave-front set of $\pi$ is 
\[
 \operatorname{WF}(\pi)
 =\bigcup_{c_{\OO}(\pi)\ne0}\clan{\OO}.
\]
We choose $\mathcal T_t$ which meets only one of the orbits 
in $ \operatorname{WF}(\pi)$ which is maximal.

\begin{lemma}\label{lem:orbit-isolation}
Let $\OO'\in  \operatorname{WF}(\pi)$ be maximal, and let $s\in\OO'$.
Fix the lattice $\Lambda$ and integer $m$ chosen in
Section~\ref{sec:positive}, and let
\[
 Q=\exp(\varpi^m\Lambda),\qquad
 L=\begin{cases}
     \varpi^m\Lambda,&p\ne2,\\
     2\varpi^m\Lambda,&p=2.
   \end{cases}
\]
Let $L^\perp$ denote the annihilator defined in
\eqref{eqn: lattice-perp}.
For all sufficiently small $t\in F^\times$, the compact open subset 
$\mathcal T_t=\Ad(Q)(t^{-2}s)+L^\perp$ of $\g$ meets no other orbit
than $\OO'$ in $\operatorname{WF}(\pi)$.
\end{lemma}

\begin{proof}
The finite union
\[
 \mathcal Z=
 \bigcup_{\substack{\OO''\ne\OO'\\c_{\OO''}(\pi)\ne0}}
      \clan{\OO''}
\]
is closed and $G$-invariant, and $s\notin\mathcal Z$.
Choose a full lattice $M_0\subset\g$ with
$s+M_0\subset\g\setminus\mathcal Z$.
For sufficiently small $t\in F^\times$, compactness of $L^\perp$ gives
\[
 t^2L^\perp\subset M_0,\qquad
 s+t^2L^\perp\subset\g\setminus\mathcal Z.
\]
The lattice $L^\perp$ is preserved by $\Ad(Q)$, as proved in
Proposition~\ref{prop:exact-transporter}. Hence
\[
 \begin{aligned}
 t^2\mathcal T_t
 =\Ad(Q)s+t^2L^\perp=\Ad(Q)(s+t^2L^\perp)
 \subset\Ad(Q)(\g\setminus\mathcal Z)
 =\g\setminus\mathcal Z.
 \end{aligned}
\]

For a nonzero rational nilpotent $x$, the neutral cocharacter
$\gamma_{h_x}$ supplied by Section~\ref{subsec:neutral-model} satisfies
\[
 \Ad(\gamma_{h_x}(t^{-1}))x=t^2x.
\]
Thus $t^2\OO''=\OO''$ for every rational nilpotent orbit.
For each $\OO''\ne\OO'$ with $c_{\OO''}(\pi)\ne0$, we conclude
\[
 \mathcal T_t\cap\OO''
 =t^{-2}(t^2\mathcal T_t\cap\OO'')=\varnothing.
\]
\end{proof}

\begin{proof}[Proof of Theorem~\ref{thm:main}]
For $\pi=0$, one has $\Wh_{\OO}(\pi)=0$ and
$\operatorname{WF}(\pi)=\varnothing$, so the equivalence holds.
Assume henceforth that $\pi\ne0$.
The forward implication is Theorem~\ref{thm:published-forward}.
Suppose $\OO\subset\operatorname{WF}(\pi)$.
If $\OO=\{0\}$, then $\Wh_0(\pi)=\Hom_{\C}(\pi,\C)\ne0$.

Now take $0\ne f\in\OO$. The definition of
$\operatorname{WF}(\pi)$ and finiteness of the analytic closure
order give an orbit $\OO'$, maximal among those occurring in
$\operatorname{WFC}(\pi)$, such that
\[
 f\in\clan{\OO'},\qquad
 c_{\OO'}(\pi)\ne0.
\]
Lemmas~\ref{lem:slice-incidence} and~\ref{lem:adapted-lattice}
give a point $s\in\OO'\cap S_f$ and an adapted lattice $\Lambda$;
the replacement of $s$ in the second lemma stays in $\OO'$.
Choose $m$ as in Section~\ref{sec:positive}, so that
$Q=\exp\Lambda_m\subset\mathcal U$.
Choose $t\in F^\times$ sufficiently small for both
Proposition~\ref{prop:positive-test} and
Lemma~\ref{lem:orbit-isolation}.

Suppose $\Wh_f(\pi)=0$. Proposition~\ref{prop:positive-test}
supplies $a_t\in C_c^\infty(Q)$ with
\[
 \begin{gathered}
 \pi(a_t)=0,\qquad \widehat{a_t^\natural}\ge0,\qquad
 \supp(\widehat{a_t^\natural})\subset\mathcal T_t,
 \end{gathered}
\]
and 
\[
 \widehat{a_t^\natural}(t^{-2}s)>0,
 \qquad t^{-2}s\in\OO'.
\]
The set $\mathcal T_t$ is open in $\g$, being a union of cosets
of the open lattice $L^\perp$.
Choose a compact open neighborhood $D\subset\mathcal T_t$ of $t^{-2}s$
in $\g$ on which $\widehat{a_t^\natural}$ is constant.
Then $D\cap\OO'$ contains $t^{-2}s$ and is open in the relative
topology of $\OO'$. Full support of $\mu_{\OO'}$ therefore gives
$\mu_{\OO'}(D\cap\OO')>0$.
Since $\widehat{a_t^\natural}\in C_c^\infty(\g)$,
\eqref{eq:orbital-integrability} applied to
$b=\widehat{a_t^\natural}$ and $\OO=\OO'$ gives
\[
 0<
 \widehat{a_t^\natural}(t^{-2}s)\,\mu_{\OO'}(D\cap\OO')
 \le
 \int_{\OO'}\widehat{a_t^\natural}(\xi)\,d\mu_{\OO'}(\xi)
 <\infty.
\]
For every $\OO''\ne\OO'$ with $c_{\OO''}(\pi)\ne0$,
Lemma~\ref{lem:orbit-isolation}, together with the support inclusion
in Proposition~\ref{prop:positive-test}, gives
\[
 \OO''\cap\supp(\widehat{a_t^\natural})
 \subset\OO''\cap\mathcal T_t=\varnothing.
\]
Thus $\widehat{a_t^\natural}$ vanishes at every point of $\OO''$, and
\[
 \int_{\OO''}\widehat{a_t^\natural}(\xi)\,d\mu_{\OO''}(\xi)=0.
\]
Since $a_t$ is supported in $\mathcal U$, the character expansion
\eqref{eq:character-expansion} therefore reduces to
\[
 0=\Tr\pi(a_t)
   =c_{\OO'}(\pi)
      \int_{\OO'}\widehat{a_t^\natural}(\xi)\,d\mu_{\OO'}(\xi)
   \ne0.
\]
The last inequality uses $c_{\OO'}(\pi)\ne0$ and the strictly
positive integral. This contradiction proves $\Wh_{\OO}(\pi)\ne0$.
\end{proof}

\bibliographystyle{amsplain}
\bibliography{biblio}

\end{document}